\documentclass[10pt]{article}
\usepackage{amsmath}
\usepackage{graphicx}
\usepackage{multirow}
\usepackage{epsfig}
\usepackage{color}
\def\disp{\displaystyle}

\def\Limsup{\mathop{{\rm Lim}\,{\rm sup}}}

\def\tto{\;{\lower 1pt \hbox{$\rightarrow$}}\kern -10pt
\hbox{\raise 2pt \hbox{$\rightarrow$}}\;}
\def\Hat{\widehat}

\def\Bar{\overline}
\def\ra{\rangle}
\def\la{\langle}
\def\ve{\varepsilon}
\def\B{I\!\!B}
\def\h{\hfill\Box}
\def\R{I\!\!R}
\def\N{I\!\!N}
\def\ox{\bar{x}}

\def\dom{\mbox{\rm dom}\,}

\def\substack#1#2{{\scriptstyle{#1}\atop\scriptstyle{#2}}}

\def\h{\hfill\triangle}
\def\dn{\downarrow}

\def\ph{\varphi}
\def\emp{\emptyset}
\def\st{\stackrel}
\def\oR{\Bar{\R}}

\def\gg{\gamma}
\def\dd{\delta}
\def\al{\alpha}

\def \N{I\!\!N}

\newcounter{lk}

\begin{document}
\begin{center}
\vspace*{0.3in} {\bf APPLICATIONS OF VARIATIONAL ANALYSIS\\ TO A
GENERALIZED HERON PROBLEM}
\\[2ex]
BORIS S. MORDUKHOVICH\footnote{Department of Mathematics, Wayne
State University, Detroit, MI 48202, USA and King Fahd University of Petroleum and Minerals, Dhahran, Saudi Arabia (email:
boris@math.wayne.edu).},
NGUYEN MAU NAM\footnote{Department of Mathematics, The University of Texas--Pan American, Edinburg, TX
78539--2999, USA (email: nguyenmn@utpa.edu).} and JUAN SALINAS JR.\footnote{Department of Mathematics, The University of
Texas--Pan American, Edinburg, TX 78539--2999, USA (email:
jsalinasn@broncs.utpa.edu).}\\[2ex]
\end{center}
\small{\bf Abstract.} This paper is a continuation of our ongoing efforts to solve a number of
geometric problems and their extensions by using advanced tools of
variational analysis and generalized differentiation. Here we
propose and study, from both qualitative and numerical viewpoints, the following optimal
location problem as well as its further extensions: on a given nonempty subset of a Banach
space, find a point such that the sum of the distances from it to
$n$ given nonempty subsets of this space is minimal. This is a
generalized version of the classical Heron problem: on a given
straight line, find a point $C$ such that the sum of the distances
from $C$ to the given points $A$ and $B$ is minimal. We show that the advanced variational techniques
allow us to completely solve optimal location problems of this type in some important settings.
\\[1ex]

\newtheorem{Theorem}{Theorem}[section]
\newtheorem{Proposition}[Theorem]{Proposition}
\newtheorem{Remark}[Theorem]{Remark}
\newtheorem{Lemma}[Theorem]{Lemma}
\newtheorem{Corollary}[Theorem]{Corollary}
\newtheorem{Definition}[Theorem]{Definition}
\newtheorem{Example}[Theorem]{Example}
\renewcommand{\theequation}{\thesection.\arabic{equation}}
\normalsize

\section{Introduction and Problem Formulation}
\setcounter{equation}{0}

In this paper we propose and largely investigate various extensions
of the Heron problem, which seem to be mathematically interesting
and important for applications. In particular, the one of this type
is to replace two given points in the classical Heron problem by
finitely many nonempty closed subsets of a Banach space and to
replace the straight line therein by another nonempty closed subset
of this space. The reader are referred to our paper \cite{mns} for partial results concerning a convex
version of this problem in the Euclidean space $\R^n$.\vspace*{0.05in}

Recall that the classical Heron problem was posted by Heron from Alexandria (10--75 AS)
in his {\em Catroptica} as follows: find a point on a straight line in the plane such that
the sum of the distances from it to two given points is minimal; see \cite{cr,h} for more discussions.
We formulate the {\em distance function version} of the
{\em generalized Heron problem} as follows:
\begin{equation}\label{distance function}
\mbox{minimize }\;D(x):=\sum_{i=1}^n d(x;\Omega_i)\;\mbox{ subject to
}\;x\in\Omega,
\end{equation}
where $\Omega$ and $\Omega_i$, $i=1,\ldots,n, n\ge 2$, are given nonempty
closed subsets of a Banach space $X$ endowed with the norm $\|\cdot\|$, and where
\begin{equation}\label{df}
d(x;Q):=\inf\big\{\|x-y\|\;\big|\;y\in Q\big\}.
\end{equation}
is the usual distance from $x\in X$ to a set $Q$. Observe that in this new formulation the
generalized Heron problem \eqref{distance function} is an extension
of the {\em generalized Fermat-Torricelli} problem proposed and
studied in \cite{mnft}. The difference is that the latter problem in
unconstrained, i.e., $\Omega=X$ in \eqref{distance function} while the
presence of the {\em geometric constraint} in the generalized Heron
version \eqref{distance function} makes it more mathematically
complicated and more realistic for applications. Among the most
natural areas of applications we mention constrained problems
arising in location science, optimal networks, wireless
communications, etc. We refer the reader to the corresponding discussions and results in \cite{mnft} and the bibliographies therein concerning unconstrained Fermat-Torricelli-Steiner-Weber versions. Needless to say that the presence of geometric (generally nonconvex) constraints in \eqref{distance function} essentially changes these versions while referring us to the original Heron geometric problem. \vspace*{0.05in}

In fact, we are able to investigate a more general version of problem
\eqref{distance function}, where the distance function \eqref{df} is
replaced by the so-called {\em minimal time function}
\begin{equation}\label{minimal time}
T^F_Q(x):=\inf\big\{t\ge 0\big|\;Q\cap(x+tF)\ne\emp\big\}
\end{equation}
with the {\em constant dynamics} $\dot x\in F\subset X$ and the {\em
target} set $Q\subset X$ in a Banach space $X$; see \cite{bmn10} and
the references therein for more discussions and results on this
class of functions important for various aspects of optimization theory and
its numerous applications.

The main problem under consideration in this paper, called below the
{\em generalized Heron problem}, is formulated as follows:
\begin{equation}\label{ft}
\mbox{minimize }\;T(x):=\sum_{i=1}^n T^F_{\Omega_i}(x)\;\mbox{ subject
to }\;x\in\Omega,
\end{equation}
where $F$ is a closed, bounded, and convex set containing the origin
as an interior point, and where $\Omega$ and $\Omega_i$ for $i=1,\ldots,n$
are nonempty closed subsets of a Banach space $X$; these are the
{\em standing assumptions} of the paper.

When $F=\B$ in \eqref{ft}, this problem reduces to the one in
\eqref{distance function}. Note that involving the minimal time
function \eqref{minimal time} into \eqref{ft} instead of the
distance function in \eqref{distance function} allows us to cover
some important location models that cannot be encompassed by
formalism \eqref{distance function}; cf.\ \cite{npr} for the case of
convex unconstrained problems of type \eqref{ft} and \cite{mnft} for
the generalized Fermat-Torricelli problem corresponding to
\eqref{ft} with $\Omega=X$.\vspace*{0.05in}

A characteristic feature of the generalized Heron problem \eqref{ft}
and its distance function specification \eqref{distance function} is
that they are {\em intrinsically nonsmooth}, since the functions \eqref{df} and
\eqref{minimal time} are nondifferentiable. These problems are generally
nonconvex while the convexity of both cost functions in \eqref{distance function} and
\eqref{ft} follows from the convexity the sets $\Omega_i$. This
makes it natural to apply advanced methods and tools of variational
analysis and generalized differentiation to study these problems. To
proceed in this direction, we largely employ the recent results from
\cite{bmn10} on generalized differentiation of the minimal time
function \eqref{minimal time} in convex and nonconvex settings as
well as comprehensive rules of generalized differential calculus. As
can be seen from the solutions below, the constraint nature of the
Heron problem and its extensions leads to new structural phenomena
in comparison with the corresponding Fermat-Torricelli
counterparts. Note that a number of the results obtained in this paper are new even for
the unconstrained setting of the generalized Fermat-Torricelli problem. \vspace*{0.05in}

The rest of the paper is organized as follows. In Section~2, we present some basic constructions and
properties from variational analysis that are widely used in the sequel. Section~3 concerns deriving necessary optimality conditions for solutions to the
generalized Heron problem in the case of arbitrary closed sets $\Omega$
and $\Omega_i$, $i=1,\ldots,n$, in \eqref{ft} and its specification
\eqref{distance function}. The results obtained are expressed in
terms of the limiting normal cone to closed sets in the sense
of Mordukhovich \cite{mor06a}. We pay a special attention to the
Hilbert space setting, which allows us to establish necessary (in some cases necessary and sufficient)
optimality conditions in the most efficient forms. Some examples are
given to illustrate applications of general results in particular
situations. In Section~4 we develop a numerical algorithm to solve some versions of the generalized Heron problem in finite dimensions while the concluding Section~5 is devoted to the implementation of this algorithm and its specifications in various settings of their own interest.\vspace*{0.05in}

Our notation is basically standard in the area of variational
analysis and generalized differentiation; see \cite{mor06a,rw}. We
recall some of them in the places they appear.

\section{Tools of Generalized Differentiation}
\setcounter{equation}{0}

This section contains basic constructions and results of the
generalized differentiation theory in variational analysis employed in what follows.
The reader can find all the proofs, discussions, and additional material in the books
\cite{bz,mor06a,mor06b,rw,s} and the references therein.

Given an extended-real-valued function $\ph\colon X\to\oR:=(-\infty,\infty]$ with $\ox$ from the domain $\dom\ph:=\{x\in X|\;\ph(x)<\infty\}$ and given $\ve\ge 0$, define first the
$\ve$-{\em subdifferential} of $\ph$ at $\ox$ by
\begin{eqnarray}\label{2.3}
\Hat\partial_\ve\ph(\ox):=\Big\{x^*\in
X^*\Big|\;\liminf_{x\to\ox}\frac{\ph(x)-\ph(\ox)-\la
x^*,x-\ox\ra}{\|x-\ox\|}\ge-\ve\Big\}.
\end{eqnarray}
For $\ve=0$ the set $\Hat\partial\ph(\ox):=\Hat\partial_0\ph(\ox)$
is known as {\em Fr\'echet/regular subdifferential} of
$\ph$ at $\ox$. It follows from definition (\ref{2.3}) that regular subgradients are described
as follows: $x^*\in \Hat\partial_\ve\ph(\ox)$ if and only if for any
$\eta>0$ there is $\gg>0$ such that
\begin{equation*}
\la x^*,x-\ox\ra\le\ph(x)-\ph(\ox)+(\ve+\eta)\|x-\ox\|\;\mbox{ whenever }\;x\in\ox+\gg\B
\end{equation*}
with $\B$ standing for the closed unit ball of the space in question. When $\ph$ is Fr\'echet differentiable at $\ox$, its regular subdifferential $\Hat\partial\ph(\ox)$ reduces to the classical gradient $\{\nabla\ph(\ox)\}$. Despite the simple definition \eqref{2.3} closely related to the classical derivative, the regular subdifferential and its $\ve$-enlargements in general do not happen to be
appropriate for applications to the generalized Heron problem under consideration due to the serious lack of calculus rules.

To get a better construction, we need to employ a certain robust limiting procedure, which lies at the heart of variational analysis.
Recall that, given a set-valued mapping $G\colon X\tto X^*$ between a Banach space $X$ and its topological dual $X^*$, the {\em sequential Painlev\'e-Kuratowski outer limit} of $G$ as $x\to\ox$ is defined by
\begin{eqnarray}\label{pk}
\begin{array}{ll}
\disp\Limsup_{x\to\ox}G(x):=\Big\{x^*\in X^*\Big|&\exists\,\mbox{
sequences } \;x_k\to\ox,\;x^*_k\st{w^*}{\to}x^*\;\mbox{ as }
\;k\to\infty\\
&\mbox{such that }\;x^*_k\in G(x_k)\;\mbox{ for all
}\;k\in\N:=\{1,2,\ldots\}\Big\},
\end{array}
\end{eqnarray}
where $w^*$ signifies the weak$^*$ topology of $X^*$. Applying the limiting operation \eqref{pk} to the set-valued mapping $(x,\ve)\tto\Hat\partial_\ve\ph(x)$ in
 \eqref{2.3} and using the notation $x\st{\ph}{\to}\ox:=x\to\ox$ with $\ph(x)\to\ph(\ox)$ give us the subgradient set
\begin{eqnarray}\label{2.4}
\partial\ph(\bar x):=\Limsup_\substack{x\xrightarrow{\ph}\bar x}{\ve\dn 0}
\Hat\partial_\ve\ph(x)
\end{eqnarray}
known as the {\em Mordukhovich/limiting subdifferential} of $\ph$ at $\ox$. We can equivalently put $\ve=0$ in \eqref{2.4} if $\ph$ is lower
semicontinuous around $\ox$ and if $X$ is {\em Asplund}, i.e., each of its separable subspaces has a separable dual; the latter is automatics, e.g., when $X$ is reflexive. Recall that $\ph$ is {\em subdifferentially regular} at $\ox$ if $\partial\ph(\ox)=\Hat\partial\ph(\ox)$.

Note that every convex function $\ph$ is subdifferentially regular at any point $\ox\in\dom\ph$ with the classical subdifferential representation
\begin{equation}\label{cs}
\partial\ph(\ox)=\big\{x^*\in X^*\;\big|\;\la x^*,
x-\ox\ra\le\ph(x)-\ph(\ox)\;\mbox{ for all }\;x\in X\big\}.
\end{equation}
However, the latter property often fails in nonconvex setting, where $\Hat\partial\ph(\ox)$ may be empty (as for $\ph(x)=-|x|$ at $\ox=0$) with a poor calculus, while the limiting subdifferential \eqref{2.4} enjoys a {\em full calculus} (at least in Asplund spaces) due to {\em variational/extremal principles} of variational analysis. We following calculus results are most useful in this paper.

\begin{Theorem}\label{sum rule} {\bf (subdifferential sum rules).} Let $\ph_i\colon X\to\oR$,
$i=1,\ldots,n$, be lower semicontinuous functions on a Banach space $X$. Suppose that all but one of them are locally Lipschitzian
around $\ox\in\cap_{i=1}^n\dom\ph_i$. Then:

{\bf (i)} We have the inclusion
\begin{equation}\label{sum1}
\partial\Big(\sum_{i=1}^n\ph_i\Big)(\ox)\subset\sum_{i=1}^n\partial\ph_i(\ox)
\end{equation}
provided that $X$ is Asplund. Furthermore, inclusion \eqref{sum1} becomes an equality if all the functions $\ph_i$ are subdifferentially regular at $\ox$.

{\bf (ii)} When all the functions $\ph_i$ are convex, the equality
\begin{equation}\label{sum2}
\partial\Big(\sum_{i=1}^n\ph_i\Big)(\ox)=\sum_{i=1}^n\partial\ph_i(\ox)
\end{equation}
holds with no Asplund space requirement.
\end{Theorem}

Note that assertion (ii) of Theorem~\ref{sum rule}, which is the classical Moreau-Rockafellar theorem, is a consequence of assertion (i) in the case of Asplund spaces; see \cite[Theorem~3.36]{mor06a}. \vspace*{0.05in}

Finally in this section, recall that the corresponding {\em normal cones} to a set $\Omega$ at $\ox\in\Omega$ can be defined via the subdifferentials \eqref{2.3} and \eqref{2.4} of the indicator function by
\begin{equation}\label{nc}
\Hat N(\ox;\Omega):=\Hat\partial\dd(\ox;\Omega)\;\mbox{ and }\;N(\ox;\Omega):=\partial\dd(\ox;\Omega),
\end{equation}
where $\dd(x;\Omega)=0$ if $x\in\Omega$ and $\dd(x;\Omega)=\infty$ otherwise.

\section{Optimality Conditions for the Generalized Heron Problem}
\setcounter{equation}{0}

The main results of this section give necessary optimality conditions for the generalized Heron problem under consideration, which occur to be necessary and sufficient for optimality in the case of convex data. To begin with, we would like make sure that problem \eqref{ft} admits an optimal solution under natural assumptions.

\begin{Proposition}\label{exi} {\bf (existence of optimal solutions to the generalized Heron problem).} The generalized Heron problem \eqref{ft} admits an optimal solution in each of the following three cases:

{\bf (i)} $X$ is a Banach space, and the constraint set $\Omega$ is compact.

{\bf (ii)} $X$ is finite-dimensional, and one of the sets $\Omega$ and $\Omega_i$ as $i=1,\ldots,n$ is bounded.

{\bf (iii)} $X$ is reflexive, the sets $\Omega$ and $\Omega_i$ as $i=1,\ldots,n$ are
convex and one of them  is bounded.
\end{Proposition}
{\bf Proof.}It follows from \cite[Proposition~2.2]{mn10} that the minimal time function (\ref{minimal time}) and hence the function $T$ in \eqref{ft} are Lipschitz continuous. Thus the conclusion in the case (i) follows from the classical Weierstrass theorem.

Consider the infimum value
\begin{equation*}
\gamma:=\inf_{x\in\Omega}T(x)<\infty
\end{equation*}
in problem (\ref{ft}) and take a minimizing sequence $\{x_k\}$ with $T(x_k)\to\gamma$ as
$k\to\infty$ and $x_k\in\Omega$ for all $k\in\N$. Now assume that $X$ is finite dimensional and $\Omega_1$ is bounded. When $k$ is sufficiently large, one has
\begin{equation*}
T^F_{\Omega_1}(x_k)\leq T(x_k)<\gamma +1.
\end{equation*}
Thus there exist $0\leq t_k<\gamma+1$, $f_k\in F$, and $w_k\in
\Omega_1$ such that
\begin{equation*}
x_k+t_kf_k=w_k.
\end{equation*}
Since both $F$ and $\Omega_1$ are bounded, $(x_k)$ is a bounded
sequence, and hence it has subsequence that converges to $\ox\in \Omega$.
Then $\ox$ is a solution of the problem under (ii). The proof in case (iii) is similar to that given in \cite[Proposition~4.1]{mns}. $\h$\vspace*{0.05in}

To proceed with deriving optimality conditions for the generalized Heron problem \eqref{ft} and its specification \eqref{distance function}, we need more notation. Define the {\em support level set}
\begin{equation*}
C^*:=\big\{x^*\in X^*\big|\;\sigma_F(-x^*)\le 1\big\}
\end{equation*}
via the {\em support function} of the constant dynamics
\begin{equation*}
\sigma_F(x^*):=\sup_{x\in F}\la x^*,x\ra,\quad x^*\in X^*.
\end{equation*}
The {\em generalized projection} to the target set $Q$ via the minimal time function \eqref{minimal time} is a set-valued mapping $\Pi^F_Q\colon X\tto X$ defined by
\begin{equation}\label{pr}
\Pi^F_Q(x):=Q\cap\big(x+T^F_Q(x)F\big),\quad x\in X.
\end{equation}
Considering further the {\em Minkowski gauge}
\begin{equation}\label{mg}
\rho_F(x):=\inf\big\{t\ge 0\big|\;x\in t F\big\},\quad x\in X,
\end{equation}
and involving the limiting normal cone from \eqref{nc}, we define the sets
\begin{eqnarray}\label{A}
A_i(x):=\left\{\begin{array}{ll}
\disp\bigcup_{\omega\in\Pi^F_{\Omega_i}(x)}\big[-\partial
\disp\rho_F(\omega-x)\cap N(\omega;\Omega_i)\big]\;\mbox{ for }\;x\notin\Omega_i,\;\Pi^F_{\Omega_i}(x)\ne\emp,\\\\
N(x;\Omega_i)\cap C^*\;\mbox{ for }\;x\in\Omega_i\;\mbox{ as }\;i=1,\ldots,n.
\end{array}\right.
\end{eqnarray}
We say that the minimal time function $T^F_{Q}(\cdot)$ is {\em well
posed} at $\ox$ if for every sequence $\{x_k\}$ converging to
$\ox$ there is a sequence $\{w_k\}$ such that $w_k\in\Pi^F_Q(x_k)$
and $\{w_k\}$ contains a convergent subsequence. The reader is
referred to \cite[Proposition~6.2]{bmn10} for a number of verifiable conditions ensuring such a well-posedness of the minimal time function.

Our first theorem establishes necessary as well as necessary and sufficient conditions for optimality in \eqref{ft} via the sets $A_i(x)$ from \eqref{A} in general infinite-dimensional settings.

\begin{Theorem}\label{necessary} {\bf (necessary and sufficient optimality conditions for the generalized Heron problem in Banach and Asplund spaces).} Given $\ox\in\Omega$, suppose in the setting of \eqref{ft} that the minimal time function $T^F_{\Omega_i}$ is well posed at $\ox$ for each $i\in\{1,\ldots,n\}$ such that $\ox\notin\Omega_i$. The following assertions hold:

{\bf (i)} Let $\ox$ be a local optimal solution to \eqref{ft}, and let $X$ be Asplund. Then we have
\begin{equation}\label{necessary1}
0\in\sum_{i=1}^n A_i(\ox)+N(\ox;\Omega),
\end{equation}
where the sets $A_i(\ox)$ are defined in \eqref{A}.

{\bf (ii)} Let $X$ be a general Banach space, and let all the sets $\Omega$ and $\Omega_i$ as $i=1,\ldots,n$ be convex. Given $\ox\in\Omega$, assume that
$\Pi^F_{\Omega_i}(\ox)\ne\emp$ for $i=1,\ldots,n$ with $\ox\notin\Omega_i$, select any $\omega\in\Pi^F_{\Omega_i}(\ox)$, and construct $A_i(\ox)$ by
\begin{eqnarray}\label{c}
A_i(\ox):=N(\bar\omega;\Omega_i)\cap\big[-\partial\rho_F(\bar\omega-\ox)\big]\;\mbox{ for }\;\ox\notin\Omega_i
\end{eqnarray}
and by the second formula in \eqref{A} otherwise. Then $\ox$ is an optimal solution to \eqref{ft} if and only if inclusion \eqref{necessary1} is satisfied.
\end{Theorem}
{\bf Proof.} Observe first that problem \eqref{ft} can be equivalently written in the form
\begin{equation}\label{ft0}
\mbox{minimize }\;T(x)+\delta(x;\Omega).
\end{equation}
It easily follows from definitions \eqref{2.3} and \eqref{2.4} of regular and limiting subgradients and their description \eqref{cs} for convex functions that the generalized Fermat rule
\begin{equation}\label{fermat}
0\in\Hat\partial f(\ox)\subset\partial f(\ox)
\end{equation}
is a necessary condition for a local minimizer $\ox$ of any function $f\colon X\to\oR$ being also sufficient for this if $f$ is convex. To justify now assertion (i), we apply \eqref{fermat} via $\partial f(\ox)$ to the cost function $f(x):=T(x)+\dd(x;\Omega)$ in \eqref{ft0} and then use the subdifferential sum rule for limiting subgradients from Theorem~\ref{sum rule}(i) in Asplund spaces by taking into account that the functions $T^F_{\Omega_i}$ are Lipschitz continuous. It follows in this way that
\begin{equation}\label{nc1}
\begin{array}{ll}
0&\in\partial\big(T+\delta(\cdot;\Omega)\big)(\ox)\subset\partial T(\ox)+N(\ox;\Omega)\\
&\subset\disp\sum_{i=1}^n\partial T^F_{\Omega_i}(\ox)+N(\ox;\Omega).
\end{array}
\end{equation}
Employing further the subdifferential formulas for the minimal time function from \cite[Theorem~3.1 and Theorem~3.2]{mnft} gives us
\begin{equation}\label{nc2}
\partial T^F_{\Omega_i}(\ox)\subset A_i(\ox),\quad i=1,\ldots,n.
\end{equation}
Substituting the latter into \eqref{nc1} justifies inclusion \eqref{necessary1} in assertion (i) of the theorem.

To justify assertion (ii), we apply Theorem~\ref{sum rule}(ii) for convex functions on Banach spaces and conclude in this way that both inclusions ``$\subset$" in \eqref{nc1} hold as equalities and provide necessary and sufficient optimality conditions for optimality of $\ox$ in \eqref{ft}. Employing finally \cite[Theorem~7.1 and 7.3]{bmn10} gives us the equalities in \eqref{nc2}, where the sets $A_i(\ox)$ are calculated by \eqref{c} when $\ox\notin\Omega_i$. This completes the proof of the theorem. $\h$\vspace*{0.05in}

It is not hard to check under our standing assumptions that the requirement $\Pi^F_{\Omega_i}(\ox)\ne\emp$ in Theorem~\ref{necessary}(ii) is automatically satisfied when the space $X$ is reflexive.\vspace*{0.05in}

The next theorem allows us to significantly simplify the calculation of the sets $A_i(\ox)$ in Theorem~\ref{necessary} for the case of Hilbert spaces and thus to ease the implementation of the optimality conditions obtained therein. Besides this, it leads us to an improvement of optimality under some additional assumptions. Namely, we can replace the limiting normal cone in \eqref{necessary1} by the smaller regular one for an arbitrary closed constraint set $\Omega$. Define the index sets
\begin{equation}\label{ind}
I(x):=\big\{i\in\{1,\ldots,n\}\;\big|\;x\in\Omega_i\big\}\;\mbox{ and }\;J(x)=\big\{i\in\{1,\ldots,n\}\;\big|\;x\notin\Omega_i\big\},\quad x\in X.
\end{equation}
We obviously have $I(x)\cup J(x)=\{1,\ldots,n\}$ and $I(x)\cap J(x)=\emp$ for all $x\in X$.

\begin{Theorem}\label{nonconvex} {\bf (improved optimality conditions in Hilbert spaces).} Consider version \eqref{distance function} of the generalized Heron problem with a Hilbert space $X$ in the assumptions of Theorem~{\rm\ref{necessary}}. The following assertions hold:

{\bf (i)} Let $\ox\in\Omega$ be a local optimal solution to \eqref{distance function}, and let $\Pi(\ox;\Omega_i)\ne\emp$ whenever $i\in J(\ox)$. Then for any $a_i(\ox)\in A_i(\ox)$ as $i\in J(\ox)$ we have
\begin{equation}\label{Frechet1}
-\sum_{i\in J(\ox)}a_i(\ox)\in\sum_{i\in I(\ox)}A_i(\ox)+N(\ox;\Omega),
\end{equation}
where each set $A_i(\ox)$ is computed by
\begin{eqnarray}\label{A5}
A_i(\ox)=\left\{\begin{array}{lr}
\dfrac{\ox-\Pi(\ox;\Omega_i)}{d(\ox;\Omega_i)}&\mbox{ for
}\;\ox\notin\Omega_i,\\\\
N(\ox;\Omega_i)\cap\B&\mbox{ for }\;\ox\in\Omega_i
\end{array}\right.
\end{eqnarray}
whenever $i=1,\ldots,n$. If in addition $I(\ox)=\emp$, then
\begin{equation}\label{Frechet2}
-\sum_{i=1}^na_i(\ox)\in\Hat N(\ox;\Omega).
\end{equation}

{\bf (ii)} If all the sets $\Omega$ and $\Omega_i$ as $i=1,\ldots,n$ are convex, then each set $A_i(\ox)$ as $i\in J(\ox)$ in \eqref{A5} is a singleton $\{a_i(\ox)\}$ and condition \eqref{Frechet1} is necessary and sufficient for the global optimality of $\ox\in\Omega$ in problem \eqref{distance function}.
\end{Theorem}
{\bf Proof.} To justify assertion (i), pick $\bar\omega_i\in\Pi(\bar x;\Omega_i)$ for all $i\in J(\ox)$ such that $a_i(\ox)=\dfrac{\ox-\bar\omega_i}{d(\ox;\Omega_i)}$ and get the relationships
\begin{equation*}
\sum_{i\in J(\ox)}\|\ox-\bar\omega_i\|+\sum_{i\in I(\ox)}d(\ox;\Omega_i)=\sum_{i=1}^n d(\ox;\Omega_i)\le\sum_{i=1}^n d(x;\Omega_i)\le\sum_{i\in J(\ox)}\|x-\bar\omega_i\|+\sum_{i\in I(\ox)} d(x;\Omega_i)
\end{equation*}
for all $x\in\Omega$ around $\ox$. This shows that $\ox$ is a local optimal solution to the problem
\begin{equation}\label{sm}
\mbox{minimize }\;p(x):=\sum_{i\in J(\ox)}\|x-\bar\omega_i\|+\sum_{i\in I(\ox)}d(x;\Omega_i)\;\mbox{ subject to }\;x\in\Omega.
\end{equation}
Since the norm function on a Hilbert space is Fr\'echet differentiable in any nonzero point, we conclude that each
$p_i(x):=\|x-\bar\omega_i\|$ as $i\in J(\ox)$ is Fr\'echet differentiable at $\ox$ with
\begin{equation*}
\nabla p_i(\ox)=\dfrac{\ox-\bar\omega_i}{\|\ox-\bar\omega_i\|}=\dfrac{\ox-\bar\omega_i}{d(\ox;\Omega_i)}=a_i(\ox).
\end{equation*}
Applying to \eqref{sm} the first inclusion in the generalized Fermat rule \eqref{fermat} and then using the subdifferential sum rules from \cite[Proposition~1.107(i)]{mor06a} for regular subgradients and from Theorem~\ref{sum rule}(i) for limiting ones, we get
\begin{align*}
0\in\Hat\partial\big[p+\delta(\cdot;\Omega)\big](\ox)&=\sum_{i\in J(\ox)}\nabla p_i(\ox)+\Hat\partial\Big[\sum_{i\in I(\ox)}d(\cdot;\Omega_i)+\delta(\cdot;\Omega)\Big](\ox)\\
&\subset\sum_{i\in J(\ox)}a_i(\ox)+\partial\Big[\sum_{i\in I(\ox)}d(\cdot;\Omega_i)+\delta(\cdot;\Omega)\Big](\ox)\\
&\subset\sum_{i\in J(\ox)}a_i(\ox)+\sum_{i\in I(\ox)}\partial d(\ox;\Omega_i)+N(\ox;\Omega)\\
&\subset\sum_{i\in J(\ox)}a_i(\ox)+\sum_{i\in I(\ox)}[N(\ox;\Omega_i)\cap\B]+N(\ox;\Omega)\\
&=\sum_{i\in J(\ox)}a_i(\ox)+\sum_{i\in I(\ox)}A_i(\ox)+N(\ox;\Omega),
\end{align*}
where the last three relationships hold since $\ox\in\Omega_i$ for each $i\in I(\ox)$. This justifies inclusion (\ref{Frechet1}). In the case of $I(\ox)=\emp$, we arrive at inclusion (\ref{Frechet2}) by the first row of the above relationships and the normal cone definition \eqref{nc}.

Assertion (ii) is justified similarly to the proof of Theorem~\ref{necessary}(ii) by using the results of assertion (i) and the well-known fact that the projection operator for a closed and convex set in a Hilbert space is single-valued. $\h$\vspace*{0.05in}

Observe that in Theorem~\ref{nonconvex}, in contrast to Theorem~\ref{necessary}, we do not impose the well-posedness requirement. In fact, under the assumptions of Theorem~\ref{nonconvex}(ii) it holds automatically; see \cite[Corollary~1.106]{mor06a}. Note also that in finite-dimensional spaces $X$ we always have the Fr\'echet differentiability of the distance function at out-of-set points with unique projections (see, e.g., \cite[Exercise~8.53]{rw}), and so we can deal in the proof of Theorem~\ref{nonconvex}(i) directly with the cost function in the generalized Heron problem \eqref{distance function}, without considering the auxiliary problem \eqref{sm}. However, in Hilbert spaces this approach requires additional and unavoidable assumptions on the projection continuity; see \cite[Corollary~3.5]{f}. In finite dimensions the projection continuity and Fr\'echet differentiability of the distance functions actually follows from the projection uniqueness, while it is not the case in Hilbert spaces as shown in \cite[Example~5.2]{f}. Observe to this end that neither uniqueness nor continuity of projections is required in Theorem~\ref{nonconvex}.\vspace*{0.05in}

On the other hand, the next result shows that for the unconstrained version of \eqref{distance function}, i.e., for the generalized Fermat-Torricelli problem \cite{mnft} with disjoint sets $\Omega_i$, the projection nonemptiness at a local optimal solution automatically implies the projection uniqueness in arbitrary Hilbert spaces.

\begin{Proposition}\label{singleton} {\bf (projection uniqueness at optimal solutions).} Let $\ox$ be a local optimal solution to problem \eqref{ft} in a Hilbert space $X$ with $\Omega=X$ and $\cap_{i=1}^n\Omega_i=\emp$. Assume that $\ox\notin\Omega_i$ as $i=1,\ldots,n$. Then the fulfillment of the condition $\Pi(\ox;\Omega_i)\ne\emp$ for all $i=1,\ldots,n$ implies that the projection set $\Pi(\ox;\Omega_i)$ is a singleton whenever $i\in\{1,\ldots,n\}$.
\end{Proposition}
{\bf Proof.} Since $I(\ox)=\emp$ for the first index set in \eqref{ind}, it follows from the proof of Theorem~\ref{nonconvex}(i) with $\Omega=X$ that for every $\omega_i\in\Pi(\ox;\Omega_i)$ as $i=1,\ldots,n$ we have the equality
\begin{equation}\label{sm1}
0=\sum_{i=1}^n\dfrac{\ox-\omega_i}{d(\ox;\Omega_i)}.
\end{equation}
Picking any $\Omega_i$, say $\Omega_1$, let us check that the set $\Pi(\ox;\Omega_1)$ is singleton. Indeed, take two projections $\bar\omega_{1,1}, \bar\omega_{1,2}\in\Pi(\ox;\Omega_1)$ and fix arbitrary projections $\bar\omega_i\in\Pi(\ox;\Omega_i)$ for $i=2,\ldots,n$. Then from \eqref{sm1} we get the relationships
\begin{equation*}
0=\dfrac{\ox-\bar\omega_{1,1}}{d(\ox;\Omega_1)}+\sum_{i=2}^n\dfrac{\bar x-\bar\omega_i}{d(\ox;\Omega_i)}=\dfrac{\ox-\bar\omega_{1,2}}{d(\ox;\Omega_1)}+\sum_{i=2}^n\dfrac{\bar x-\bar\omega_i}{d(\ox;\Omega_i)},
\end{equation*}
which imply that $\bar\omega_{1,1}=\bar\omega_{1,2}$ and thus complete the proof of the proposition. $\h$\vspace*{0.05in}

Observe that if $\ox$ belongs to one of the sets $\Omega_i$ as $i=1,\ldots,n$, the conclusion of Proposition~\ref{singleton} does not generally hold even in finite dimensions as it is demonstrated by the following example.

\begin{Example}\label{proj} {\bf (nonuniqueness of projections at solution points).}  {\rm Let $X=\R^2$ in the setting of Proposition~\ref{singleton}, let $\Omega_1$ be the unit circle of $\R^2$, and let $\Omega_2=\{(0,0)\}$. Then $\ox=\{(0,0)\}$ is a solution of the Fermat-Torricelli problem generated by $\Omega_1$ and $\Omega_2$, but the projection $\Pi(\ox;\Omega_1)$ is the whole unit circle. It is also clear that any point inside of the unit circle other than $(0,0)$ is also a solution to this problem, and $\Pi(\ox;\Omega_i)$ is a singleton for both $i=1,2$, which is consistent with the result of Proposition~\ref{singleton}.}
\end{Example}

The observation made in Proposition~\ref{singleton} allows us to improve the optimality conditions obtained in \cite[Corollary~4.1]{mnft} for the generalized Fermat-Torricelli problem.

\begin{Corollary}\label{necessary conditions} {\bf (improved optimality conditions for the generalized Fermat-Torricelli problem with three nonconvex
sets in Hilbert spaces).} Let $n=3$ in the framework of Theorem~{\rm\ref{nonconvex}}, where $\Omega_1,\Omega_2$, and $\Omega_3$ are pairwisely
disjoint subsets of $X$ and $\Omega=X$. The following alternative holds for a local optimal solution $\ox\in X$ with the sets $A_i(\ox)$ defined by
\eqref{A5}:

{\bf (i)} The point $\ox$ belongs to one of the sets $\Omega_i$, say
$\Omega_1$. Then for any $a_i\in A_i(\ox)$ as $i=2,3$ we have the relationships
\begin{equation*}
\big\la a_2,a_3\ra\le-1/2\;\mbox{ and }\;-a_2-a_3\in \Hat N(\ox;\Omega_1).
\end{equation*}

{\bf (ii)} The point $\ox$ does not belong to all the three sets
$\Omega_1$, $\Omega_2$, and $\Omega_3$. Then $A_i(\ox)=\{a_i\}$ for all
$i=1,2,3$ and we have
\begin{equation*}
\la a_i,a_j\ra=-1/2\;\mbox{ for }\;i\ne j\;\mbox{ as
}\;i,j\in\big\{1,2,3\big\}.
\end{equation*}
Conversely, suppose that the sets $\Omega_i$, $i=1,2,3$, are convex and that $\ox$ satisfies either {\rm(i)} or {\rm (ii)}. Then it is a global optimal solution to the problem under consideration.
\end{Corollary}
{\bf Proof.} In case (i) for any $a_i\in A_i(\ox)$ as $i=2,3$ take $\bar\omega_i\in\Pi(\ox;\Omega_i)$ such that
\begin{equation*}
a_i=\dfrac{\ox-\bar\omega_i}{d(\ox;\Omega_i)},\quad i=2,3.
\end{equation*}
Since $\ox\in\Omega_1$, we have the relationships
\begin{equation*}
\|\bar x-\bar\omega_1\|+\|\bar x-\omega_2\|=\sum_{i=1}^3 d(\ox;\Omega_i)\le\sum_{i=1}^3 d(x;\Omega_i)\le d(x;\Omega_1)+\|x-\bar\omega_2\|+\|x-\bar\omega_3\|
\end{equation*}
whenever $x$ is near $\ox$. Thus $\ox$ is a local optimal solution to the problem
\begin{equation}\label{sm2}
\mbox{minimize }\;q(x):=d(x;\Omega_1)+\|x-\bar\omega_2\|+\|x-\bar\omega_3\|.
\end{equation}
Employing the generalized Fermat rule in \eqref{sm2} and then the aforementioned sum rule for regular subgradients gives us by using the well-known formula for the regular subdifferential of the distance function (see, e.g., \cite[Corollary~1.96]{mor06a}) that
\begin{equation*}
0\in\Hat\partial q(\ox)=\Hat\partial d(\ox;\Omega_1)+a_2+a_3=\Hat N(\ox;\Omega_1)\cap\B+a_2+a_3.
\end{equation*}
The latter implies therefore that
\begin{equation*}
-a_2-a_3\in\Hat N(\ox;\Omega_1)\;\mbox{ with }\;\|a_2+a_3\|\le 1.
\end{equation*}
The rest of the proof follows the lines of that in \cite[Corollary~4.1]{mnft}. Assertion (ii) and the converse statement are derived similarly from Proposition~\ref{singleton} and the proof of \cite[Corollary~4.1]{mnft} by the same procedure, which thus allows us to fully justify the corollary. $\h$\vspace*{0.05in}

From now on in this section we concentrate on the distance function version \eqref{distance function} of the generalized Heron problem while paying the main attention to deriving efficient forms of optimality conditions for \eqref{distance function} under additional structural assumptions on the constraint set $\Omega$.  In what follows in this section we impose the {\em nonintersection condition}
\begin{equation}\label{noint}
\Omega\cap\Omega_i=\emp\;\mbox{ for all }\;i=1,\ldots,n
\end{equation}
on the sets $\Omega$ and $\Omega_i$ in \eqref{distance function}, which is specific for the (constrained) generalized Heron problem. In this case we obviously have $I(\ox)=\emp$ for the first index set in \eqref{ind} whenever $\ox\in\Omega$, and so the sets $A_i(\ox)$ are calculated by
\begin{equation}\label{a1}
A_i(\ox)=\dfrac{\ox-\Pi(\ox;\Omega_i)}{d(\ox;\Omega_i)},\quad i=1,\ldots,n,
\end{equation}
in the Hilbert space setting under consideration.

To proceed, for any nonzero vectors $u,v\in X$ define the quantity
\begin{equation*}
\cos(u, v):=\dfrac{\la u,v\ra}{\|u\|\cdot\|v|}
\end{equation*}
and, given a linear subspace $L$ of $X$, recall that
\begin{equation*}
L^\perp:=\big\{x^*\in X\big|\;\la x^*,v\ra=0\;\mbox{ for all }\;v\in
L\big\}.
\end{equation*}
We say that $\Omega\subset X$ has a {\em tangent space} $L=L(\ox)$ at $\ox$ if $L^\perp=\Hat N(\ox;\Omega)$. Note that for any affine subspace $\Omega\subset X$ parallel to a linear subspace $L$ the tangent space to $\Omega$ at every $\ox\in\Omega$ is $L$.

Next we derive verifiable necessary and sufficient conditions for optimal solutions to \eqref{distance function} in Hilbert spaces provided that the constraint set admits a tangent space at the reference point.

\begin{Proposition}\label{cos} {\bf (optimality conditions for the case of constraint sets with tangent spaces).} Consider the generalized Heron problem \eqref{distance function} under condition \eqref{noint} in Hilbert spaces. The following assertions hold:

{\bf (i)} Let $\ox\in\Omega$ be a local optimal solution to \eqref{distance function}, let $A_i(\ox)$ be computed in \eqref{a1} where \\$\Pi(\ox;\Omega_i)\neq\emptyset$ for $i=1,\ldots,n$, and let $\Omega$ admit a tangent space $L(\ox)$ at $\ox$. Then for any
$a_i(\ox)\in A_i(\ox)$, one has
\begin{equation}\label{cos rep}
\sum_{i=1}^n\cos\big(a_i(\ox),v\big)=0\;\mbox{ for every }\;v\in L(\ox)\setminus\{0\}.
\end{equation}

{\bf (ii)} Let all the sets $\Omega_i$, $i=1,\ldots,n$, be convex. Then $A_i(\ox)=\{a_i(\ox)\}$ and condition \eqref{cos rep} with the tangent space $L(\ox)$ for $\Omega$ is necessary and sufficient for the global optimality of $\ox$ in \eqref{distance function}.
\end{Proposition}
{\bf Proof.} To justify (i), observe by the assumptions made and the definition of the tangent space  $L(\ox)$ to $\Omega$ at $\ox$ that
\begin{equation*}
\Hat N(\ox;\Omega)=L^\perp=\big\{v\in X\big|\;\la v,x\ra=0\;\mbox{ for all }\;x\in L(\ox)\big\}.
\end{equation*}
By Theorem \ref{nonconvex} for any $a_i(\ox)\in A_i(\ox)$, one has
\begin{equation*}
0\in\sum_{i=1}^n a_i(\ox)+L^\perp(\ox),
\end{equation*}
which implies in turn that
\begin{equation*}
\Big\la\sum_{i=1}^n a_i(\ox),v\Big\ra=0\;\mbox{ for all }\;v\in L(\ox).
\end{equation*}
Since $\ox\notin\Omega_i$ by \eqref{noint}, we have due to \eqref{a1} that $\|a_i(\ox)\|=1$ for $i=1,\ldots,n$, and
hence
\begin{equation*}
\sum_{i=1}^n \dfrac{\la a_i(\ox),v\ra}{\|a_i(\ox)\|\cdot\|v\|}=0\;\mbox{ whenever }\;v\in L(\ox)\setminus\{0\}.
\end{equation*}
Thus we arrive at the the necessary optimality condition \eqref{cos rep}.

To justify (ii), observe that the implication ``$\Longrightarrow$" follows directly from assertion (i) of the theorem, since
the sets
$A_i(\ox)$ are singletons for $i=1,\ldots,n$  in this case. The oppositive implication ``$\Longleftarrow$" follows from
Theorem~\ref{nonconvex}(ii) by taking into account the special structure of the normal cone $\Hat N(\ox;\Omega)=L^\perp(\ox)$.
This completes the proof of the proposition. $\h$\vspace*{0.05in}

We have the following specification of optimality conditions in Proposition~\ref{cos} when the tangent space therein is finitely generated.

\begin{Corollary}\label{finite dimensional} {\bf (optimality conditions for the case of finitely generated tangent spaces).} Let $L(\ox)={\rm span}\{v_1,\ldots,v_s\}$ with $v_j\ne 0$ as $j=1,\ldots,s$ in the setting of Proposition~{\rm\ref{cos}}. Then condition
\eqref{cos rep} in all of its conclusions is equivalent to
\begin{equation}\label{cos rep1}
\sum_{i=1}^n\cos (a_i,v_j)=0\;\mbox{ for all }\;j=1,\ldots,s.
\end{equation}
\end{Corollary}
{\bf Proof.} We obviously have that (\ref{cos rep})$\Longrightarrow$(\ref{cos
rep1}). To justify the converse implication, set $a:=\sum_{i=1}^n a_i$ and observe by
$v_j\ne 0$ as $j=1,\ldots,s$ and $\|a_i\|=1$ as
$i=1,\ldots,n$ that (\ref{cos rep1}) yields $\la a,
v_j\ra=0$ for all $j=1,\ldots,s$. Picking further an arbitrary vector $v\in L(\ox)\setminus\{0\}$,
we arrive at the representation
\begin{equation*}
v=\sum_{j=1}^s\lambda_j v_j
\end{equation*}
with some $\lambda_j\in\R$. It gives by linearity that $\la a,v\ra=\sum_{j=1}^n\lambda_j\la
a,v_j\ra =0$, which yields (\ref{cos rep}) and completes the proof of the proposition. $\h$\vspace*{0.05in}

The next result concerns the generalized Heron problem for two nonconvex sets in Hilbert spaces with a one-dimensional structure of the regular normal cone to the constraint.

\begin{Proposition}\label{two set} {\bf (necessary conditions for the generalized Heron problem with two nonconvex sets in Hilbert spaces).} Consider problem \eqref{distance function} for two sets $(n=2)$ in Hilbert spaces under the nonintersection condition \eqref{noint}. Let $\ox\in\Omega$ be a local optimal solution to \eqref{distance function} such that
$\Hat N(\ox;\Omega)={\rm span}\{v\}$ with some $v\ne 0$ and that $\Pi(\ox;\Omega_i)\ne\emp$ for $i=1,2$. Then for any $a_i(\ox)\in A_i(\ox)$ as $i=1,2$ we have the conditions:
\begin{equation}\label{v*}
\mbox{either }\;a_1(\ox)+a_2(\ox)=0\;\mbox{ or }\;\cos\big(a_1(\ox),
v\big)=\cos\big(a_2(\ox),v\big).
\end{equation}
\end{Proposition}
{\bf Proof.} It follows from Theorem~\ref{nonconvex}(i) in this setting that
\begin{equation}\label{v0}
-a_1(\ox)-a_2(\ox)\in\Hat N(\ox;\Omega)\;\mbox{ for any }\;a_i(\ox)\in A_i(\ox),\;i=1,2.
\end{equation}
Denoting for simplicity $a_i:=a_i(\ox)$ as $i=1,2$ and taking into account the assumed structure of the regular normal cone to $\Omega$, we get that \eqref{v0} is equivalent to the following:
\begin{equation*}
\mbox{either }\;a_1+a_2=0\;\mbox{ or }\;
a_1+a_2=\lambda v \mbox{ with some }\;\lambda\ne 0.
\end{equation*}
Let us show that the latter condition implies that $\cos(a_1,v)=\cos(a_2,v)$.
Indeed, in this case we have $\|a_1\|=\|a_1\|=1$, which gives by the Euclidean norm on $X$ that
\begin{align*}
\lambda^2 \|v\|^2=\|a_1+a_2\|^2=\|a_1\|^2+\|a_2\|^2+2\la a_1,
a_2\ra=2+2\la a_1,a_2\ra.
\end{align*}
This implies in turn the relationships
\begin{align*}
\la a_1,\lambda v\ra&=\la\lambda v-a_2,\lambda v\ra\\
&=\lambda^2\|v\|^2-\lambda\la a_2,v\ra\\
&=2+2\la a_1,a_2\ra-\lambda\la a_2,v\ra\\
&=2\la a_2,a_2\ra+2\la a_1,a_2\ra-\lambda\la a_2,v\ra\\
&=2\la a_2+a_1,a_2\ra-\lambda\la a_2,v\ra\\
&=2\la\lambda v,a_2\ra-\lambda\la a_2,v\ra=\la a_2,\lambda
v\ra,
\end{align*}
which yield that $\la a_1,v\ra=\la a_2,v\ra$ since $\lambda\ne 0$. By  taking into account that $\|a_1\|=\|a_2\|=1$ and
$v\ne 0$, we conclude that $\cos(a_1,v)=\cos(a_2,v)$ and thus complete the proof. $\h$\vspace*{0.05in}

Observe that sufficient optimality conditions in the form of Proposition~\ref{two set} do not hold even in convex settings. The next result provides slightly modified conditions, which are sufficient for optimality in the case of the convex generalized Heron problem on the plane.

\begin{Proposition}\label{suff2} {\bf (characterizing optimal solutions for the generalized Heron problem with two convex sets).} Let the sets $\Omega_1$ and $\Omega_2$ be convex in the setting of Proposition~{\rm\ref{two set}}, and let $a_i:=a_i(\ox)$ as $i=1,2$. Then the modification
\begin{equation}\label{normalvector1}
\mbox{either }\;a_1+a_2=0\;\mbox{ or }\;\big[a_1\ne a_2\;\mbox{ and
}\;\cos(a_1,v)=\cos(a_2,v)\big],
\end{equation}
of the necessary condition \eqref{v*} is sufficient for the global optimality of $\ox\in\Omega$ in \eqref{distance function} when $X=\R^2$.
\end{Proposition}
{\bf Proof.} To justify the sufficiency of conditions \eqref{normalvector1} for the optimality of $\ox$ in \eqref{distance function}, we need to show---by taking into account Theorem~\ref{nonconvex}(ii) and the assumed structure of the regular normal cone to $\Omega$---that the
relationships in (\ref{normalvector1}) imply the fulfillment of
\begin{equation}\label{nv1}
-a_1-a_2\in\Hat N(\ox;\Omega)={\rm{span}}\{v\}.
\end{equation}
When $-a_1-a_2=0$, inclusion \eqref{nv1} is obviously satisfied. Consider
the alternative in (\ref{normalvector1}) when $a_1\ne a_2$ and
$\cos(a_1,v)=\cos(a_2,v)$. Since we are in $\R^2$, represent
$a_1=(x_1,y_1)$, $a_2=(x_2,y_2)$, and $v=(x,y)$ with two real
coordinates. Then the equality $\cos(a_1,v)=\cos(a_2,v)$ can be
written as
\begin{equation}\label{dotprod}
x_1x+y_1y=x_2x+y_2y,\;\mbox{ i.e., }\;(x_1-x_2)x=(y_2-y_1)y.
\end{equation}
Since $v\ne 0$, assume without loss of generality that $y\ne 0$. By the equivalence
\begin{equation*}
\|a_1\|^2=||a_2||^2\Longleftrightarrow x_1^2+y_1^2=x_2^2+y_2^2
\end{equation*}
we have the equality $(x_1-x_2)(x_1+x_2)=(y_2-y_1)(y_2+y_1)$, which
implies by \eqref{dotprod} that
\begin{equation}\label{dp}
y(x_1-x_2)(x_1+x_2)=x(x_1-x_2)(y_2+y_1).
\end{equation}
Note that $x_1\ne x_2$, since otherwise we have from \eqref{dotprod}
that $y_1=y_2$, which contradicts the condition $a_1\ne a_2$ in
\eqref{normalvector1}. Dividing  both sides of \eqref{dp} by
$x_1-x_2$, we get
\begin{equation*}
y(x_1+x_2)=x(y_2+y_1),
\end{equation*}
which implies in turn that
\begin{equation*}
y(a_1+a_2)=y(x_1+x_2,y_1+y_2)=\big(x(y_1+y_2),
y(y_1+y_2)\big)=(y_1+y_2)v.
\end{equation*}
In this way we arrive at the representation
\begin{equation*}
a_1+a_2=\dfrac{y_1+y_2}{y}v
\end{equation*}
showing that inclusion \eqref{nv1} is satisfied. This ensures the
optimality of $\ox$ in \eqref{distance function} and thus completes
the proof of the proposition. $\h$\vspace*{0.05in}

We conclude this section by a simple example showing how the results obtained allow us to completely solve a direct generalization of the classical Heron problem in $\R^2$, where the constraint straight line is replaced by a convex set.

\begin{Example}\label{e1} {\bf (complete solution of a a convex set extension of the Heron problem on the plane).}
{\rm Consider problem (\ref{distance function}), where $\Omega$
is the epigraph  of the nonsmooth convex function $y=|x|$ in $\R^2$, and where $\Omega_1$ and
$\Omega_2$ are two points $(x_1,y_1)$ and $(x_2,y_2)$ that do not lie on $\Omega$. This problem admits optimal solutions due to Proposition~\ref{exi}(ii). To solve it, we are going to employ appropriate necessary optimality conditions obtained above. Observe first that
the normal cone to $\Omega$ at $(0,0)$ is given by
$$
N\big((0,0);\Omega\big)=\big\{(x,y)\in\R^2\big|\;y\le-|x|\big\}
$$
while the classical normals at other points of $\Omega$ are calculated trivially. Using this, we can easily check that if the points $(x_1,y_1)$ and $(x_2,y_2)$ belong to the region
$$
\big\{(x,y)\in\R^2\big|\;y\le-|x|\big\},
$$
then the origin $\ox=(0,0)$ is the only point that satisfies the necessary optimality condition from Theorem~\ref{nonconvex}(i) written now as:
\begin{equation*}
-a_1-a_2\in N(\ox;\Omega)\;\mbox{ with }\;a_i=\dfrac{(x_i,y_i)}{\|(x_i,y_i)\|}\;\mbox{ as }\;i=1,2.
\end{equation*}
If the points $(x_1,y_1)$ and $(x_2,y_2)$ belong to another region
\begin{equation*}
\big\{(x,y)\in\R^2\big|\;x>|y|\big\},
\end{equation*}
then the problem also has a unique optimal solution constructed by
connecting the reflection point of $(x_1,y_1)$ through the line $y=x$ and
$(x_2,y_2)$.}
\end{Example}

\section{Subgradient Algorithm in the Generalized Heron Problem}
\setcounter{equation}{0}

In this section we develop a subgradient algorithm for the numerical solution of the generalized Heron problem \eqref{ft} for finitely many convex sets and convex constraints in
 the finite-dimensional Euclidean space $\R^m$. These are our standing assumptions for the rest of the paper. Recall that $\Pi(x;\Omega)$ denotes the (unique) Euclidean projection of $x$ to $\Omega$ while $\Pi^F_{\Omega_i}(x)$ stands for the generalized/minimal time projection \eqref{pr} of this point to the target sets $\Omega_i$ in \eqref{ft}. Here is the algorithm whose various implementations are presented in the next section.

\begin{Theorem}\label{subgradient method2} {\bf (subgradient algorithm for the generalized Heron problem).}
Let $S\ne\emp$ be the set of optimal solutions to problem \eqref{ft}.
Picking a sequence $\{\al_k\}_{k\in\N}$ of positive numbers and a
starting point $x_1\in\Omega$, consider the algorithm
\begin{equation}\label{al}
x_{k+1}=\Pi\Big(x_k-\al_k\sum_{i=1}^n q_{ik};\Omega\Big),\quad k=1,2,\ldots,
\end{equation}
with an arbitrary choice of vectors
\begin{equation}\label{a1a}
q_{ik}\in-\partial\rho_F(\omega_{ik}-x_k)\cap
N(\omega_{ik};\Omega_i)\;\mbox{ for some
}\;\omega_{ik}\in\Pi^F_{\Omega_i}(x_k)\;\mbox{ if }\;x_k\notin\Omega_i
\end{equation}
via the Minkowski gauge \eqref{mg} and with $q_{ik}:=0$ otherwise. Assume that
\begin{equation}\label{a2}
\sum_{k=1}^\infty\alpha_k=\infty\;\mbox{ and }\;
\ell^2:=\sum_{k=1}^\infty\alpha_k^2<\infty.
\end{equation}
Then the iterative sequence $\{x_k\}$ in \eqref{al} converges to an
optimal solution of problem \eqref{ft} and the numerical value sequence
\begin{equation}\label{Vk}
V_k:=\min\big\{T(x_j)\big|\;j=1,\ldots,k\big\}
\end{equation}
converges to the optimal value $\Hat V$ in this problem.
Furthermore, we have the estimate
\begin{align*}
V_k-\Hat
V\le\dfrac{d(x_1;S)^2+L^2\sum_{i=1}^k\alpha_i^2}{2\sum_{i=1}^k\alpha_i},
\end{align*}
where $0\le L<\infty$ is a Lipschitz constant of the function
$T(\cdot)$ from \eqref{ft} on $R^m$.
\end{Theorem}
{\bf Proof.} We know that the value function $T(\cdot)$ in
\eqref{ft} is convex and globally Lipschitzian on $\R^m$. Employing \cite[Theorems~7.1 and 7.3]{bmn10}, the convex
subdifferential of the minimal time functions \eqref{minimal time} at any point $x_k$ is computed by
\begin{eqnarray}\label{al-sub}
\partial T^F_{\Omega_i}(x_k)=\left\{\begin{array}{ll}
N(x_k;\Omega_i)\cap\big\{v\in X\big|\;\sigma_F(-v)\le 1\big\}&\mbox{if
}\;x_k\in\Omega_i,\\\\
N(\omega_{ik};\Omega_i)\cap\big[-\partial\rho_F(\omega_{ik}-x_k)\big]&\mbox{
if }\;x_k\notin\Omega_i,
\end{array}\right.
\end{eqnarray}
where $\omega_{ik}\in\Pi^F_{\Omega_i}(x_k)$ is an arbitrary generalized
projection vector for $i\in\{1,\ldots,n\}$ and $k\in\N$. Recalling now
the subgradient algorithm for minimizing the convex function
$T(\cdot)$ in \eqref{ft} subject to $x\in\Omega$, we construct the iteration sequence by
\begin{equation}\label{al-sub1}
x_{k+1}=\Pi\Big(x_k-\al_k v_k;\Omega\Big)\;\mbox{ with
}\;v_k\in\partial T(x_k),\quad k=1,2,\ldots.
\end{equation}
It follows from the convex subdifferential sum rule of Theorem~\ref{sum rule}(ii) that
\begin{equation*}
v_k=\sum_{i=1}^nq_{ik}\;\mbox{ with }\;q_{ik}\in\partial
T^F_{\Omega_i}(x_k)
\end{equation*}
for the subgradients $v_k$ in \eqref{al-sub1}. Substituting the
latter into \eqref{al-sub1} gives us algorithm \eqref{al} with
$q_{ik}$ satisfying \eqref{a1a}. Then all the conclusions of the theorem are derived from the so-called
``square summable but not summable case" of the subgradient method for constrained convex functions under the conditions in
\eqref{a2}; see \cite{bert,boyd} for more details. $\h$\vspace*{0.05in}

 In the case of $F=\B$, the closed unit ball in $\R^m$, we are able to provide a more explicit algorithm to solve the distance function version \eqref{distance function} of the generalized Heron problem with now uniquely defined vectors $q_{ik}$ in \eqref{al}.

\begin{Corollary}\label{subgradient method3} {\bf (explicit subgradient algorithm for the distance version of the generalized Heron problem).} Consider the distance function specification \eqref{distance function} of the generalized Heron problem under the assumptions of Theorem~{\rm\ref{subgradient method2}}. Then all the conclusions of this theorem hold with $q_{ik}$ in \eqref{al} calculated by
\begin{equation}\label{a1b}
q_{ik}=\left\{\begin{array}{ll} 0 &\mbox{if }\;x_k\in\Omega_i,\\\\
\dfrac{x_k-\Pi(x_k;\Omega_i)}{d(x_k;\Omega_i)} &\mbox{if }\;x_k\notin\Omega_i.
\end{array}\right.
\end{equation}
\end{Corollary}
{\bf Proof}. As follows from the proof of Theorem~\ref{nonconvex}, in the case of problem \eqref{distance function} the vectors $q_{ik}$ from \eqref{a1a} are uniquely determined and reduce to \eqref{a1b}. $\h$\vspace*{0.05in}

The next corollary specifies algorithm \eqref{al} in the case of balls for the distance function version \eqref{distance function} of the generalized Heron problem.

\begin{Corollary}\label{toballs}{\bf (subgradient algorithm in the case of multidimensional balls).}
Consider problem \eqref{distance function} with $\Omega_i=B(c_i,r_i)\subset\R^m$ as
$i=1,\ldots,n$. Then the quantities $q_{ik}$ in Theorem~{\rm\ref{subgradient method2}} are uniquely calculated by
\begin{equation}\label{a3}
q_{ik}=\left\{\begin{array}{ll} 0 &\mbox{if }\;\|x_k-c_i\|\le
r_i,\\\\
\disp\frac{x_k-c_i}{\|x_k-c_i\|} &\mbox{if }\;\|x_k-c_i\|>r_i
\end{array}\right.
\end{equation}
and the corresponding values $V_k$ are evaluated by formula
\eqref{Vk} with
\begin{equation}\label{a3a}
T(x_j)=\sum_{i=1,\;x_j\notin\Omega_i}^n\Big(\|x_j-c_i\|-r_i\Big).
\end{equation}
\end{Corollary}
{\bf Proof.} Formula \eqref{a3} directly follows from \eqref{a1b} due to the projection representation
\begin{equation*}
\Pi(x_k;\Omega_i)=c_i+r_i\dfrac{x_k-c_i}{\|x_k-c_i\|}
\end{equation*}
in the case under consideration. It is easy to see furthermore that the value function in \eqref{Vk} reduces to \eqref{a3a} in this case.
$\h$

\section{Implementation of the Subgradient Algorithm}

The final section of the paper is devoted to implementations of the subgradient algorithm from Theorem~\ref{subgradient method2} and its specifications to solve the generalized Heron problem in a number of underlying examples of their own interest. Let us start with a two-dimensional problem involving a ball constraint in the setting of Corollary~\ref{toballs}.
{\small\begin{figure}[h]
\begin{minipage}{2in}
 \includegraphics[width=3.5in]{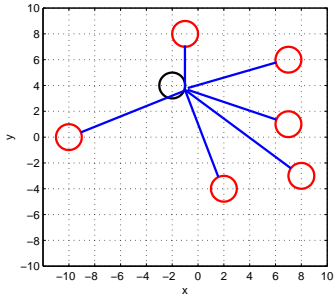}\\
\end{minipage}
~\hfill~
\begin{minipage}[t]{0.52\textwidth}
\begin{tabular}{|c|c|c|}
\hline
\multicolumn{3}{|c|}{MATLAB RESULTS} \\
\hline
$k$ & $x_k$ & $V_k$ \\
\hline
1 &       (-1,4)             & 44.58483 \\
10 &      (-1.07737,3.61433) & 44.36969 \\
100 &     (-1.07779,3.61332) & 44.36969 \\
1000 &    (-1.07779,3.61331) & 44.36969 \\
10,000 &  (-1.07779,3.61331) & 44.36969 \\
\hline
\end{tabular}
\end{minipage}
\vspace{-30pt}
\caption{A Generalized Heron Problem for Balls with a Ball
Constraint.}
\end{figure}}

\begin{Example}\label{ex4.1} {\bf (two-dimensional Heron problem for balls with ball constraints).}
{\rm Consider the generalized Heron problem \eqref{distance function} for balls in $\R^2$ subject
to a given ball constraint. Let $c_i=(a_i,b_i)$ and $r_i$ as
$i=1,\ldots,n$ be the centers and the radii of the balls $\Omega_i$
under consideration, and let $c=(x_0,y_0)$ and $r$ be the center
and radius for the given ball constraint $\Omega$. The subgradient
algorithm is given by \eqref{al}, where the projection $P(x,y):=\Pi((x,y);\Omega)$ is computed by
\begin{equation*}
P(x,y)=(v_x+x_0,v_y+y_0)\;\mbox{ with }\;v_x=\disp\frac{r(x-x_0)}{\sqrt{(x-x_0)^2+(y-y_0)^2}},\;v_y=\disp\frac{r(y-y_0)}{\sqrt{(x-x_0)^2+(y-y_0)^2}},
\end{equation*}
and where the quantities $q_{ik}$ and $V_k$ are calculated in Corollary~\ref{toballs}.

To specify the calculations, take the ball constraint $\Omega$ with center $(-2,4)$ and radius
$1$. The sets $\Omega_i,\;i=1,\ldots,6$, are the balls with centers
$(-10,0)$, $(-1,8)$, $(2,-4)$, $(7,6)$, $(7,1)$, and $(8,-3)$ and with the same
radius $r=1$. The MATLAB calculations performed by algorithm \eqref{al} with the sequence
$\al_k=1/k$ satisfying \eqref{a2} and the starting point $x_1=(-1,4)$ are presented in Figure~1. Observe that the numerical results indicate points on the ball constraint
with the optimal solution $\ox\approx(-1.07779,3.61331)$ and the
optimal value $\Hat V\approx44.36969$.}
\end{Example}

The next example concerns the generalized Heron problem with square constraints.

{\small\begin{figure}[h]
\begin{minipage}{2in}
 \includegraphics[width=4.5in]{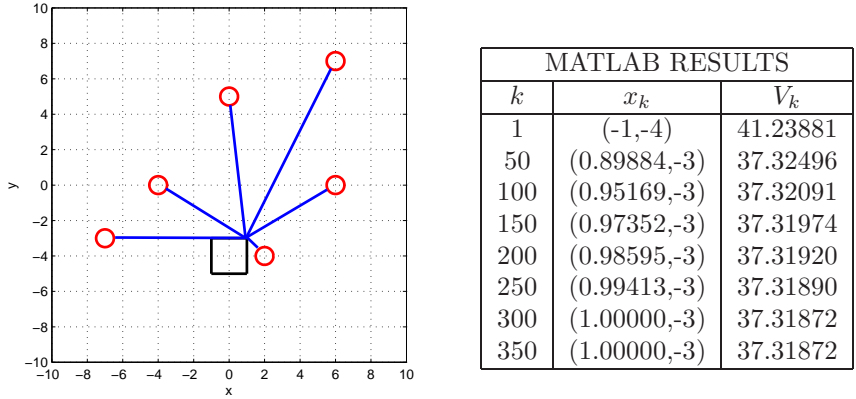}\\
\end{minipage}
~\hfill~
\begin{minipage}[t]{0.52\textwidth}
\begin{tabular}{|c|c|c|}
\hline
\multicolumn{3}{|c|}{MATLAB RESULTS} \\
\hline
$k$ & $x_k$ & $V_k$ \\
\hline
1 &   (-1,-4)      & 41.23881 \\
50 &  (0.89884,-3) & 37.32496 \\
100 & (0.95169,-3) & 37.32091 \\
150 & (0.97352,-3) & 37.31974 \\
200 & (0.98595,-3) & 37.31920 \\
250 & (0.99413,-3) & 37.31890 \\
300 & (1.00000,-3) & 37.31872 \\
350 & (1.00000,-3) & 37.31872 \\
\hline
\end{tabular}
\end{minipage}
\vspace{-20pt}
\caption{A Generalized Heron Problem for Balls with a Square
Constraint.}
\end{figure}}
\begin{Example}\label{ex4.2} {\bf (generalized Heron problem with square constraints).}
{\rm Consider the implementation of algorithm \eqref{al} for problem \eqref{distance function} using a MATLAB program with
the square constraint $\Omega$ of center $(a,b)=(0,-4)$ and short radius $r=1$
and  with the balls $\Omega_i$ as $i=1,\ldots,6$ centered at (-7,-3), (0,5),
(-4,0), (2,-4), (6,0), and (6,7) with the same radius $0.5$. Note that the projection $P(x,y)=\Pi((x,y);\Omega)$ is calculated by
{\small\begin{eqnarray*}
P(x,y)=\left\{\begin{array}{ll}
(a+r,b+r) &\mbox{if }\; x-a>r,\;y-b>r,\\\\
(x,b+r) &\mbox{if }\;|x-a|\le r,\;y-b>r,\\\\
(a-r,b+r) &\mbox{if }\;x-a<-r,\;y-b>r,\\\\
(a-r,y) &\mbox{if }\;\;x-a<-r,\; |y-b|\le r,\\\\
(a-r,b-r) &\mbox{if } x-a<-r,\;y-b<-r,\\\\
(x,b-r) &\mbox{if }\;|x-a|\le r,\;y-b<-r,\\\\
(a+r,b-r) &\mbox{if }\;x-a>r,\;y-b<-r,\\\\
(a+r,b)   &\mbox{if }\;x-a>r,\;|y-b|\le r,\\\\
(x,y) &\mbox{if }\;(x,y)\in \Omega.
\end{array}\right.
\end{eqnarray*}}
The quantities $q_{ik}$ and $V_k$ are given by Corollary~\ref{toballs}. In Figure~2 we
present the results of calculations performed by the subgradient algorithm \eqref{al} for the sequence
$\al_k=1/k$ and the starting point $x_1=(-1,-4)$. Observe that the computed optimal solution is $\ox\approx(1.00000,-3.00000)$
and the optimal value is $\Hat V\approx37.31872$.}
\end{Example}

Prior to the calculations in two next examples concerning the generalized Heron problem \eqref{distance function} for squares in $\R^2$ we formulate a specification of Theorem~\ref{subgradient method2} in a general setting of such a type. Recall that a square in $\R^2$ is of
\emph{right position} if the sides of this square are parallel to the $x$-axis and the $y$-axis, respectively.

\begin{Corollary}\label{tosquares} {\bf (subgradient algorithm for the generalized Heron problem squares targets).} Consider problem \eqref{distance function} in
$\R^2$, where each target set $\Omega_i$ is a square of right position with center $c_i=(a_i,b_i)$ and short radius $r_i$ as
$i=1,\ldots,n$, and where the constraint $\Omega$ is an arbitrary closed and convex set. Denote the vertices of the $i^{\rm th}$ square by
$v_{1i}=(a_i+r_i,b_i+r_i),\;v_{2i}=(a_i-r_i,b_i+r_i),\;v_{3i}=(a_i-r_i,b_i-r_i),\;v_{4i}=(a_i+r_i,b_i-r_i)$,
and let $x_k=(x_{1k}, x_{2k})$. Then the quantities $q_{ik}$ in Theorem~{\rm\ref{subgradient method2}} are computed by
{\small\begin{equation*}
q_{ik}=\left\{\begin{array}{ll}
0 &\mbox{if }\;|x_{1k}-a_i|\le r_i \mbox{ and }\;|x_{2k}-b_i|\le r_i,\\\\
\disp\frac{x_k-v_{1i}}{\|x_k-v_{1i}\|} &\mbox{if }\;x_{1k}-a_i>r_i\;\mbox{ and }\;x_{2k}-b_i>r_i,\\\\
\disp\frac{x_k-v_{2i}}{\|x_k-v_{2i}\|} &\mbox{if }\;x_{1k}-a_i<-r_i\;\mbox{ and }\;x_{2k}-b_i>r_i,\\\\
\disp\frac{x_k-v_{3i}}{\|x_k-v_{3i}\|} &\mbox{if }\;x_{1k}-a_i<-r_i\;\mbox{ and }\;x_{2k}-b_i<-r_i,\\\\
\disp\frac{x_k-v_{4i}}{\|x_k-v_{4i}\|} &\mbox{if }\;x_{1k}-a_i>r_i\;\mbox{ and }\;x_{2k}-b_i<-r_i,\\\\
(0,1) &\mbox{if }\;|x_{1k}-a_i|\le r_i\;\mbox{ and }\;x_{2k}-b_i>r_i,\\\\
(0,-1) &\mbox{if }\;|x_{1k}-a_i|\le r_i\;\mbox{ and }\;x_{2k}-b_i< -r_i,\\\\
(1,0) &\mbox{if }\; x_{1k}-a_i> r_i\;\mbox{ and }\;|x_{2k}-b_i|\le r_i,\\\\
(-1,0) &\mbox{if }\; x_{1k}-a_i <-r_i\;\mbox{ and }\;|x_{2k}-b_i|\le r_i\\\\
\end{array}\right.
\end{equation*}}
for all $i=1,\ldots,n$ and $k\in N$ with the corresponding quantities $V_k$ defined by \eqref{Vk}.
\end{Corollary}
{\bf Proof.} This statement follows from Corollary~\ref{subgradient method3} by a direct calculation of the projection
from an out-of-set point to each square $\Omega_i$ in formula \eqref{a1b}. $\h$\vspace*{0.05in}

Now we present the results of MATLAB calculations in the case of straight line constraints in the setting of Corollary~\ref{tosquares}.
{\small\begin{figure}[h]
\begin{minipage}{2in}
\includegraphics[width=4.0in]{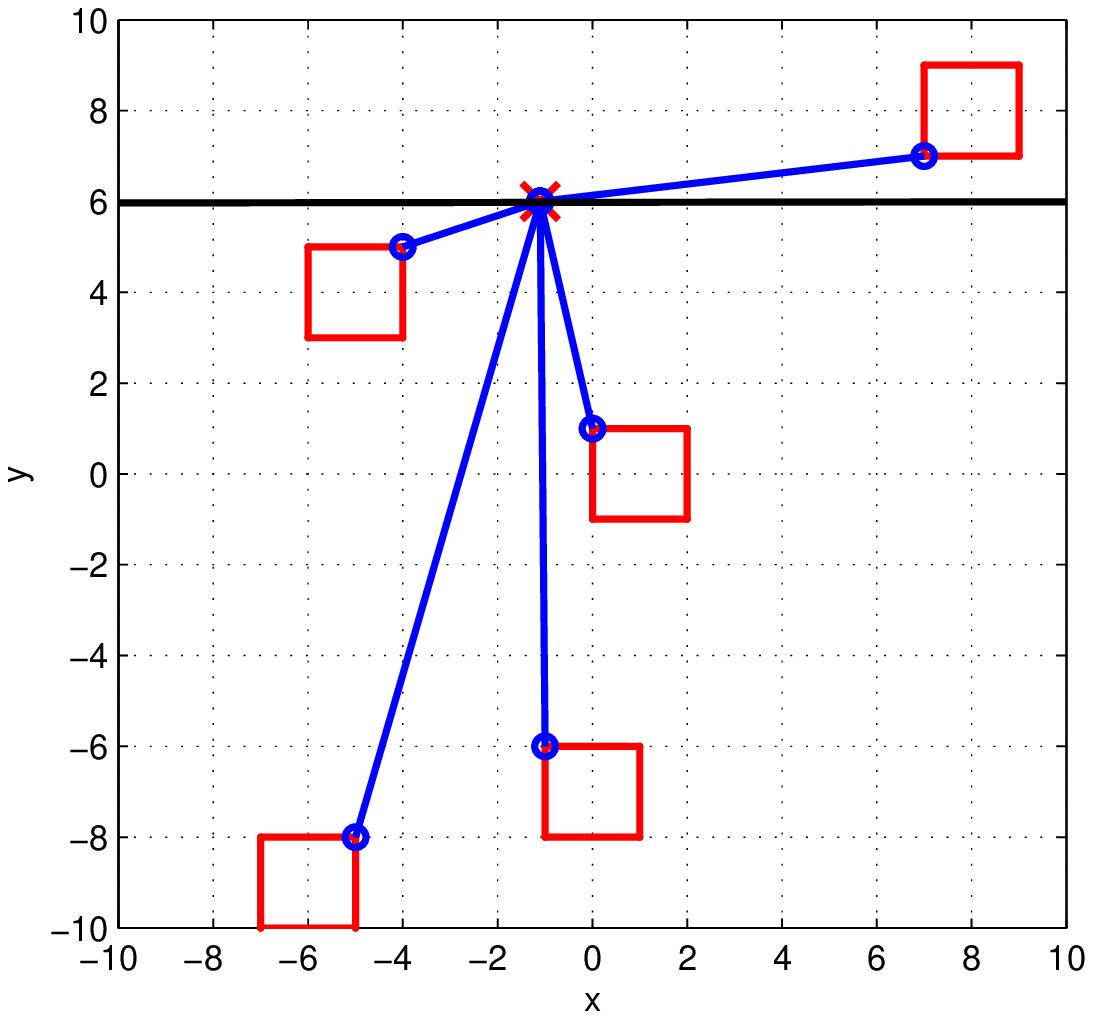}\\
\end{minipage}
~\hfill~
\begin{minipage}[t]{0.52\textwidth}
\begin{tabular}{|c|c|c|}
\hline
\multicolumn{3}{|c|}{MATLAB RESULTS} \\
\hline
$k$ & $x_k$ & $V_k$ \\
\hline
1          & (-1,6)      & 42.8838 \\
100        & (-1.0826,6) & 42.8821 \\
1000       & (-1.0896,6) & 42.8821 \\
100,000    & (-1.0938,6) & 42.8821 \\
1,000,000  & (-1.0944,6) & 42.8821 \\
5,000,000  & (-1.0946,6) & 42.8821 \\
10,000,000 & (-1.0946,6) & 42.8821 \\
\hline
\end{tabular}
\end{minipage}
\vspace{-20pt}
\caption{A Generalized Heron Problem for Squares with a Line
Constraint.}
\end{figure}}
\begin{Example}\label{ex4.3} {\bf (generalized Heron problem for squares with line constraints).}
{\rm Consider the generalized Heron problem (\ref{distance function}) for
squares of right position in $\R^2$ subject to a straight line constraint $\Omega$. Let $c_i=(a_i,b_i)$ and $r_i$ as
$i=1,\ldots,n$ be the centers and short radius of the squares $\Omega_i$ under consideration. Denote by $v_{1i}=(a_i+r_i,b_i+r_i),\;v_{2i}=(a_i-r_i,b_i+r_i),\;v_{3i}=(a_i-r_i,b_i-r_i),\;v_{4i}=(a_i+r_i,b_i-r_i)$ the vertices of the $i^{\rm th}$ square, and let
$v=[s,h]$ and $p=(x_0,y_0)$, be the direction and point vectors
of the given line $\Omega$. Then the projection $P(x,y)=\Pi((x,y);\Omega)$ in the the subgradient algorithm \eqref{al} is
calculated by
\begin{equation*}
P(x,y)=(x_0+st,y_0+ht)\;\mbox{ and }\;t=\disp\frac{s(x-x_0)+h(y-y_0)}{s^2+h^2}
\end{equation*}
while the quantities $q_{ik}$ and $V_k$ for all $i=1,\ldots,n$ and $k\in
N$ are given by Corollary~\ref{tosquares}.

In Figure~3 we present the results of calculations by algorithm \eqref{al} with $\al_k=1/k$ and the starting point $x_1=(-1,6)$ for the case above with the line constraint defined by $v=[1,0]$ and $p=(1,6)$ and the squares $\Omega_i$ as $i=1,\ldots,5$ centered at $(-6,-9)$, $(-5,4)$, $(0,-7)$, $(1,0)$, and $(8,8)$ with
the same short radius $r$=1. Observe that the calculated optimal solution is $\ox\approx(-1.0946,6)$ and the
optimal value is $\Hat V\approx42.8821$.}
\end{Example}

The next example concerns the generalized Heron problem \eqref{distance function} for squares in right position with a ball constraint on the plane.
{\small\begin{figure}[h]
\begin{minipage}{2in}
\includegraphics[width=4.8in]{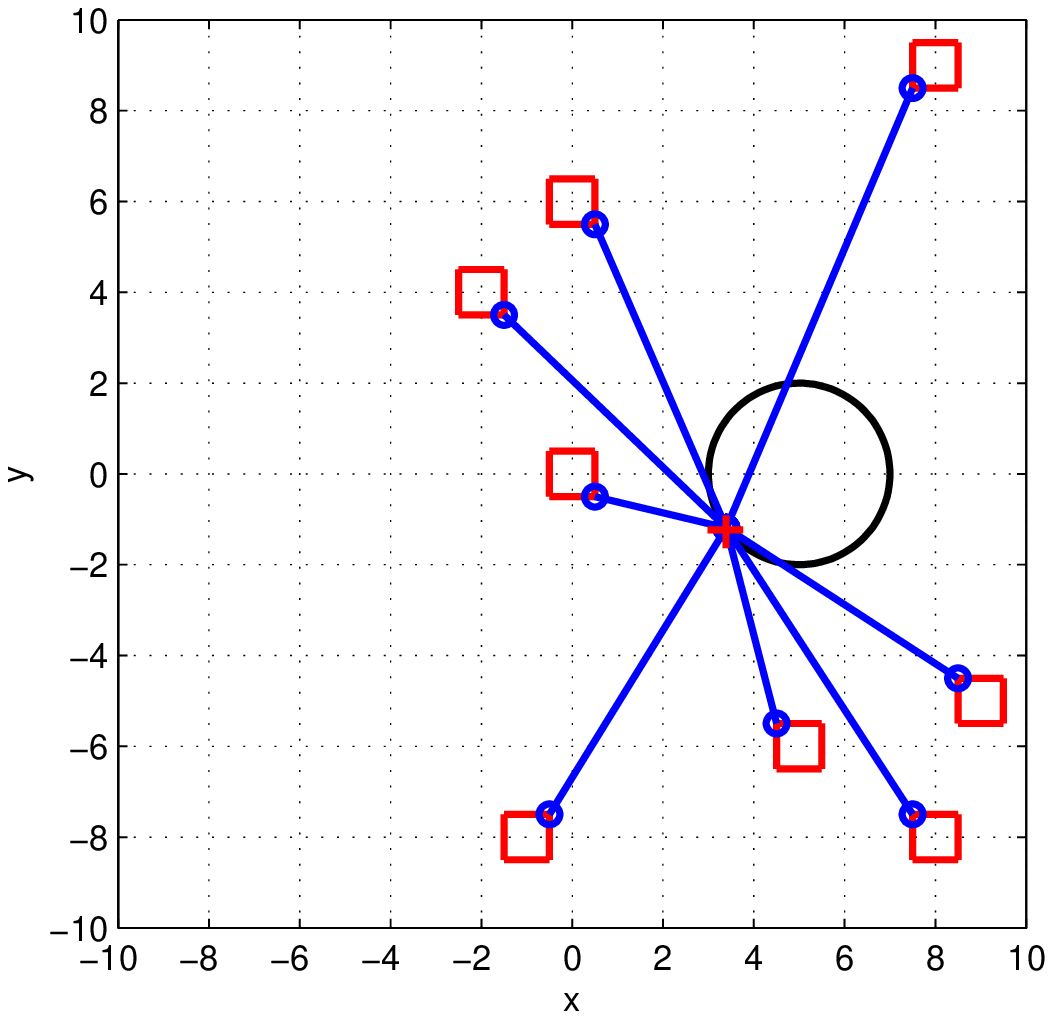}\\
\end{minipage}
~\hfill~
\begin{minipage}[t]{0.52\textwidth}
\begin{tabular}{|c|c|c|}
\hline
\multicolumn{3}{|c|}{MATLAB RESULTS} \\
\hline
$k$ & $x_k$ & $V_k$ \\
\hline
1         & (5,-2)             & 54.41891 \\
10        & (3.51379,-1.33835) & 53.05740 \\
100       & (3.41230,-1.21623) & 53.04403 \\
1000      & (3.39607,-1.19475) & 53.04364 \\
100,000   & (3.39279,-1.19033) & 53.04363 \\
600,000   & (3.39271,-1.19022) & 53.04363 \\
1,000,000 & (3.39271,-1.19021) & 53.04363 \\
1,200,000 & (3.39270,-1.19021) & 53.04363 \\
1,400,000 & (3.39270,-1.19021) & 53.04363 \\
\hline
\end{tabular}
\end{minipage}
\vspace{-20pt}
\caption{A Generalized Heron Problem for Squares with a Ball
Constraint.}
\end{figure}}
\begin{Example}\label{ex4.4} {\bf (generalized Heron problem for squares with ball constraints).} {\rm By taking into account the previous formulas for algorithm \eqref{al}, we provide the following calculations concerning the generalized Heron problem \eqref{distance function} with the ball constraint $\Omega$ centered at $(5,0)$ and radius $2$ and
the squares $\Omega_i,\;i=1,\ldots,8$ of right position with the centers $(-2,4)$, $(-1,-8)$, $(0,0)$, $(0,6)$, $(5,-6)$, $(8,-8)$, $(8,9)$, and $(9,-5)$ and the same short radius $r=0.5$. Figure~4 presents the results of calculations for algorithm \eqref{al} with the sequence $\al_k=1/k$ and the starting point $x_1= (5,-2)$. Observe that the obtained numerical results give us the optimal solution
$\ox\approx(3.39270,-1.19021)$ and the optimal value $\Hat V\approx53.04363$.}
\end{Example}

Now let us illustrate applications of the subgradient algorithm from
Theorem~\ref{subgradient method2} to solving the generalized Heron
problem \eqref{ft} formulated via the minimal time function with
dynamics sets $F$ different from the ball. First we consider the dynamics $F$ described by the
\emph{closed unit diamond}
\begin{equation}\label{dia}
F:=\big\{(x_1,x_2)\in\R^2\big|\;|x_1|+|x_2|\le 1\big\}.
\end{equation}
In this case the corresponding Minkowski gauge \eqref{mg} is given by the formula
\begin{equation}\label{mk1}
\rho_F(x_1,x_2)=|x_1|+|x_2|.
\end{equation}

The following proposition provides an explicit calculation of a subgradient of the minimal
time function \eqref{minimal time} generated by the diamond dynamics \eqref{dia} and a square target in $\R^2$. We further use this calculation in implementing algorithm \eqref{al} with the corresponding selection of $q_{ik}$ in \eqref{a1a}.

\begin{Proposition}\label{mk2} {\bf (subgradients of the minimal time function with diamond dynamics).}  Let
$F$ be the closed unit diamond in $\R^2$, and let $\Omega$ be the square
of right position centered at $c=(a,b)$ with short radius
$r>0$. Then we can calculate a subgradient $v(\ox_1,\ox_2)\in\partial
T^F_\Omega(\ox_1,\ox_2)$ of the minimal time function $T^F_\Omega$ at $(\ox_1,\ox_2)\notin\Omega$ by
{\small\begin{eqnarray}\label{v}
v(\ox_1,\ox_2)=\left\{\begin{array}{ll}
(1,0)&\mbox{if }\;|\ox_2-b|\le r,\;\ox_1>a+r,\\\\
(-1,0)&\mbox{if }\;|\ox_2-b|\le r,\;\ox_1<a-r,\\\\
(0,1)&\mbox{if }\;|\ox_1-a|\le r,\;\ox_2>b+r,\\\\
(0,-1)&\mbox{if }\;|\ox_1-a|\le r,\;\ox_2<b-r,\\\\
(1,1)&\mbox{if }\;\ox_1>a+r,\;\ox_2>b+r,\\\\
(-1,1)&\mbox{if }\;\ox_1<a-r,\;\ox_2>b+r,\\\\
(-1,-1)&\mbox{if }\;\ox_1<a-r,\;\ox_2<b-r,\\\\
(1,-1)&\mbox{if }\;\ox_1>a+r,\;\ox_2<b-r,\\\\
0 &\mbox{if }\;(\ox_1,\ox_2)\in\Omega.
\end{array}\right.
\end{eqnarray}}
\end{Proposition}
{\bf Proof}. By \cite[Theorem~7.3]{bmn10} we have the relationship
\begin{equation}\label{mk3}
\partial T^F_\Omega(\ox)=N(\bar\omega;\Omega)\cap\big[-\partial\rho_F(\bar\omega-\ox)\big]\;\mbox{ for  any }\;\bar\omega\in\Pi^F_\Omega(\ox)
\end{equation}
between the subdifferentials of the minimal time function at $\ox\notin\Omega$ and the corresponding Minkowski gauge. In the setting under consideration it is easy to find the minimal time projection $\Pi^F_\Omega(\ox_1,\ox_2)$ of a given vector $(\ox_1,\ox_2)\in\R^2$ to the square $\Omega$. Furthermore, the convex subdifferential of \eqref{mk1} at $(x_1,x_2)$ is computed by
{\small\begin{eqnarray*}
\partial\rho_F(\ox_1,\ox_2)=\left\{\begin{array}{ll}
[-1,1]\times[-1,1]&\mbox{if }\;
(\ox_1,\ox_2)=(0,0),\\\\
\disp[-1,1]\times\{1\}&\mbox{if }\;\ox_1=0,\;\ox_2>0,\\\\
\disp[-1,1]\times\{-1\}&\mbox{if }\;\ox_1=0,\;\ox_2<0,\\\\
\disp\{1\}\times[-1,1]&\mbox{if }\;\ox_1>0,\;\ox_2=0,\\\\
\disp\{-1\}\times[-1,1]&\mbox{if }\;\ox_1<0,\;\ox_2=0,\\\\
\disp\{1\}\times\{1\}&\mbox{if }\;\ox_1>0,\;\ox_2>0,\\\\
\disp\{1\}\times\{-1\}&\mbox{if }\;\ox_1>0,\;\ox_2<0,\\\\
\disp\{-1\}\times\{1\}&\mbox{if }\;\ox_1<0,\;\ox_2>0,\\\\
\disp\{-1\}\times\{-1\}&\mbox{if }\;\ox_1<0,\;\ox_2<0.\\\\
\end{array}\right.
\end{eqnarray*}}
The rest of the proof is a direct verification that the vector $v(\ox_1,\ox_2)$ from \eqref{v} belongs to the set on the right-hand side of \eqref{mk3} and hence to $\partial T^F_\Omega(\ox_1,\ox_2)$. $\h$\vspace*{0.05in}

Proposition~\ref{mk2} and the previous considerations lead us to the following realization of the subgradient algorithm \eqref{al}.

\begin{Corollary}\label{dia1} {\bf (subgradient algorithm for finitely many squares and diamond dynamics in the generalized Heron problem).}
Consider problem \eqref{ft} generated by the diamond dynamics \eqref{dia} and $n$ squares $\Omega_i$ of right position in $\R^2$. Let $c_i=(a_i,b_i)$ and $r_i$ as $i=1,\ldots,n$ be the centers and the short radii of the squares under consideration, and let
$v_{1i}=(a_i+r_i,b_i+r_i)$, $v_{2i}=(a_i-r_i,b_i+r_i)$,
$v_{3i}=(a_i-r_i,b_i-r_i)$, and $v_{4i}=(a_i+r_i,b_i-r_i)$ be the vertices of the $i^{\rm th}$ square. Denoting
$x_k=(x_{1k}, x_{2k})$ in algorithm \eqref{al}, we compute the quantities $q_{ik}$ as follows:
{\small\begin{eqnarray*}
q_{ik}=\left\{\begin{array}{ll} 0&\mbox{ if }\;|x_{1k}-a_i|\le
r_i\;\mbox{ and }\;|x_{2k}-b_i|\le
r_i,\\\\
(1,1)&\mbox{ if }\;x_{1k}-a_i>r_i\;
\mbox{ and }\;x_{2k}-b_i>r_i,\\\\
(-1,1)&\mbox{ if
}\;x_{1k}-a_i<-r_i\;\mbox{ and }\;x_{2k}-b_i>r_i,\\\\
(-1,-1)&\mbox{ if }\;x_{1k}-a_i<-r_i\;
\mbox{ and }\;x_{2k}-b_i<-r_i,\\\\
(1,-1)&\mbox{ if }\;x_{1k}-a_i>r_i\;
\mbox{ and }\;x_{2k}-b_i<-r_i,\\\\
(0,1)& \mbox{ if }\;|x_{1k}-a_i|\le r_i\;\mbox{ and
}\;x_{2k}-b_i>r_i,\\\\
(0,-1)& \mbox{ if }\;|x_{1k}-a_i|\le r_i\;\mbox{ and }
\;x_{2k}-b_i<-r_i,\\\\
(1,0)& \mbox{ if } \;x_{1k}-a_i>r_i\;\mbox{ and
}\;|x_{2k}-b_i|\le
r_i,\\\\
(-1,0)& \mbox{ if }\;x_{1k}-a_i<-r_i\;\mbox{ and
}\;|x_{2k}-b_i|\le r_i
\end{array}\right.
\end{eqnarray*}}
for all $i\in\{1,\ldots,n\}$ and $k\in\N$.
\end{Corollary}
{\bf Proof}. It follows from Proposition~\ref{mk2}, comparison between the right-hand side of \eqref{a1a} and formula \eqref{mk3}, and the square calculations of Corollary~\ref{tosquares}. $\h$\vspace*{0.05in}

Now we implement the results of Corollary~\ref{dia1} to solve the generalized Heron problem of the above type with ball constraints.
{\small\begin{figure}[h]
\vspace{-20pt}
\begin{minipage}{2in}
   \vspace{10pt}
   \hspace{-5pt}\includegraphics[width=4.2in]{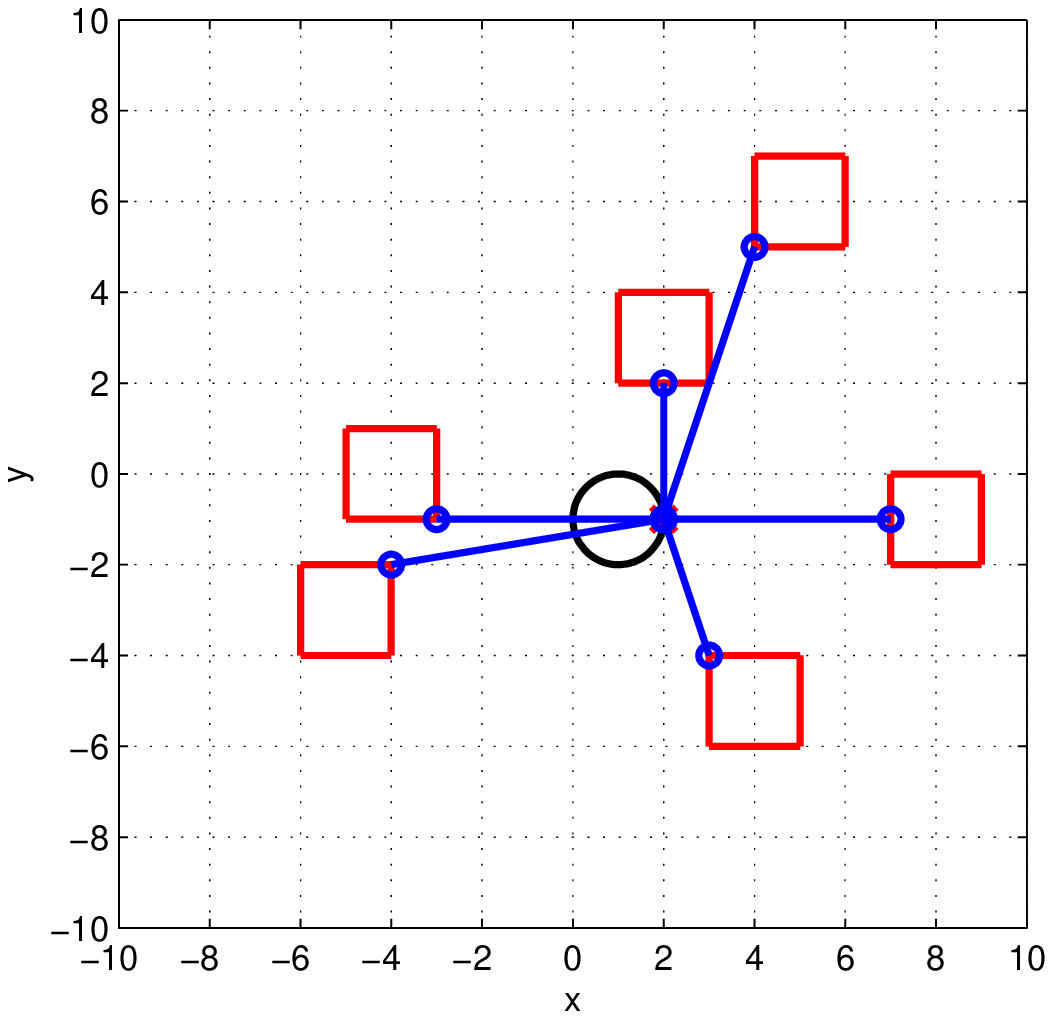}\\
\end{minipage}
~\hfill~
\begin{minipage}[t]{0.52\textwidth}
\begin{tabular}{|c|c|c|}
\hline
\multicolumn{3}{|c|}{MATLAB RESULT} \\
\hline
$k$ & $x_k$ & $V_k$ \\
\hline
1       & (1,-2)             & 34 \\
10      & (1.98703,-0.83947) & 32.01297 \\
100     & (1.99987,-0.98385) & 32.00013 \\
1,000   & (2.00000,-0.99838) & 32.00000 \\
10,000  & (2.00000,-0.99984) & 32.00000 \\
50,000  & (2.00000,-0.99997) & 32.00000 \\
100,000 & (2.00000,-0.99998) & 32.00000 \\
150,000 & (2.00000,-0.99999) & 32.00000 \\
200,000 & (2.00000,-0.99999) & 32.00000 \\
\hline
\end{tabular}
\end{minipage}
\vspace{-20pt}
\caption{A Generalized Heron Problem for Squares with a Ball
Constraint with Respect to ``Sum" Distances.}
\end{figure}}
\begin{Example}\label{ex4.5} {\bf (generalized Heron problems with diamond dynamics for squares and ball constraints).}
{\rm Consider problem \eqref{ft} with the diamond dynamics \eqref{dia} for squares $\Omega_i$ as $i=1,\ldots,6$ of right position in $\R^2$ with the centers at (-5,-3), (-4,0),
(2,3), (4,-5), (5,6), and (8,-1)  and the same short radius 1 subject to the ball constraint $\Omega$ centered at $(-1,1)$ and radius 1. The results of calculations by the subgradient algorithm \eqref{al} with $\al_k=1/k$ and the starting point $x_1=(1,-2)$ are presented in Figure~5. Observe that the obtained optimal solution is the point $\ox\approx(2.00000,-0.99999)$ on the ball constraint with the optimal value $\Hat V\approx32.00000$.}
\end{Example}
The following example is a modification of the previous one for the case of square constraints.
{\small\begin{figure}[h]
\vspace{-10pt}
\begin{minipage}{2in}
   \vspace{15pt}
   \includegraphics[width=4.2in]{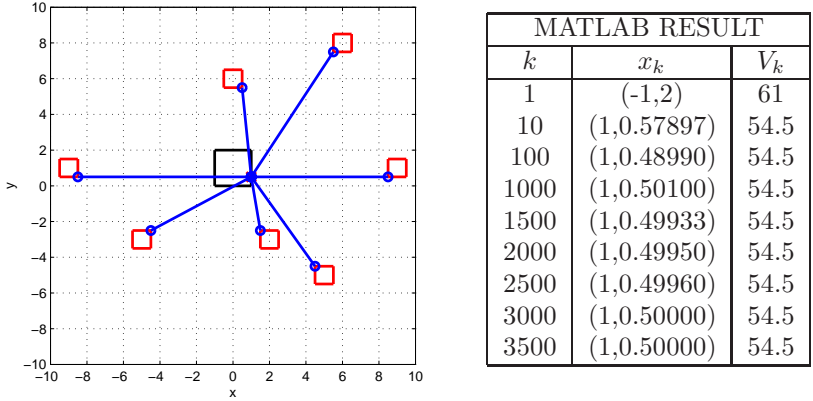}\\
\end{minipage}
~\hfill~
\begin{minipage}[t]{0.52\textwidth}
\begin{tabular}{|c|c|c|}
\hline
\multicolumn{3}{|c|}{MATLAB RESULT} \\
\hline
$k$ & $x_k$ & $V_k$ \\
\hline
1     & (-1,2)      & 61 \\
10    & (1,0.57897) & 54.5 \\
100   & (1,0.48990) & 54.5 \\
1000  & (1,0.50100) & 54.5 \\
1500  & (1,0.49933) & 54.5 \\
2000  & (1,0.49950) & 54.5 \\
2500  & (1,0.49960) & 54.5 \\
3000  & (1,0.50000) & 54.5 \\
3500  & (1,0.50000) & 54.5 \\
\hline
\end{tabular}
\end{minipage}
\vspace{-20pt}
\caption{A Generalized Heron Problem for Squares with a Square
Constraint with Respect to ``Sum" Distances.}
\end{figure}}
\begin{Example}\label{ex4.6} {\bf (generalized Heron problems for squares with diamond dynamics and square constraints).}
{\rm Consider the generalized Heron problem \eqref{ft} with the diamond dynamics \eqref{dia} for the squares $\Omega_i\in\R^2$ as $i=1,\ldots,7$ of right position centered at $(-5,-3)$, $(-9,1)$, $(0,6)$, $(2,-3)$, $(6,8)$, $(5,-5)$, and $(9,1)$ with the same short radius 1 subject to the square constraint $\Omega$ of right position centered at (0,1) with the short radius 0.5. The calculations presented in Figure~6 are performed for the sequence $\al_k=1/k$ in \eqref{al} and the
starting point $x_1=(-1,2)$. The obtained optimal solution is the point  $\ox\approx(1,0.50000)$  on the square and the
optimal value is $\Hat V\approx 54.50000$.}
\end{Example}

Next we consider the generalized Heron problem \eqref{ft} with the {\em square dynamics} $F=[-1,1]\times [-1,1]$ on the plane. The corresponding Minkowski gauge is now given by
\begin{equation*}
\rho_F(x_1,x_2)=\max\big\{|x_1|,|x_2|\big\}.
\end{equation*}
First we calculate a subgradient $v(\ox_1,\ox_2)\in\partial T^F_\Omega(\ox_1,\ox_2)$ of the cost function in \eqref{ft} at any $(\ox_1,\ox_2)$, which is further used for a specification of algorithm \eqref{al} in this setting.

\begin{Proposition}\label{mk2a} {\bf (subgradients of minimal
time functions with square dynamics and square targets).} Let
$F=[-1,1]\times[-1,1]$, and let $\Omega$ be the square of right position
in $\R^2$ centered at $c=(a,b)$ with short radius $r>0$. Then a
subgradient $v(\ox_1,\ox_2)\in\partial T^F_\Omega(\ox_1,\ox_2)$
 of the minimal time function $T^F_\Omega$
at $(\ox_1,\ox_2)$ is computed by
{\small\begin{eqnarray}\label{v1}
v(\ox_1,\ox_2)=\left\{\begin{array}{ll}
(1,0)&\mbox{if }\;|\ox_2-b|\le\ox_1-a,\;\ox_1>a+r,\\\\
(-1,0)&\mbox{if }\;|\ox_2-b|\le a-\ox_1,\;\ox_1<a-r,\\\\
(0,1)&\mbox{if }\;|\ox_1-a|\le\ox_2-b,\;\ox_2>b+r,\\\\
(0,-1)&\mbox{if }\;|\ox_1-a|\le b-\ox_2,\;\ox_2<b-r,\\\\
0 &\mbox{if }\;(\ox_1,\ox_2)\in\Omega.
\end{array}\right.
\end{eqnarray}}
\end{Proposition}
{\bf Proof.} It is given in \cite[Proposition~5.1]{mnft}. $\h$\vspace*{0.05in}

As a consequence of the proposition above, we calculate the quantities $q_{ik}$ in algorithm \eqref{al} for the corresponding version of the generalized Heron problem.

\begin{Corollary}\label{mk2b} {\bf (subgradient algorithm for the generalized Heron problem with square dynamics).} Consider problem \eqref{ft} for the square dynamics
$F=[-1,1]\times[-1,1]$ and the square targets $\Omega_i$ as $i=1,\ldots,n$ of right position in $\R^2$. Denote by $c_i=(a_i,b_i)$ and $r_i$  the
centers and the short radii of the squares $\Omega_i$ under consideration, and let
the vertices of the $i^{\rm th}$ square be
$v_{1i}=(a_i+r_i,b_i+r_i)$, $v_{2i}=(a_i-r_i,b_i+r_i)$,
$v_{3i}=(a_i-r_i,b_i-r_i)$, and $v_{4i}=(a_i+r_i,b_i-r_i)$. Then the quantities $q_{ik}$ in algorithm \eqref{al} of Theorem~{\rm\ref{subgradient method2}}
in this setting along the iterative sequence
$x_k=(x_{1k}, x_{2k})$ are calculated for all $i\in\{1,\ldots,n\}$ and $k\in\N$ by
{\small\begin{eqnarray*}
q_{ik}=\left\{\begin{array}{ll} (1,0)&\mbox{if }\;|x_{2k}-b_i|\le
x_{1k}-a_i\;\mbox{ and }\;x_{1k}>a_i+r_i,\\\\
(-1,0)&\mbox{if }\;|x_{2k}-b_i|\le a_i-x_{1k}\;\mbox{ and }\;
x_{1k}<a_i-r_i,\\\\
(0,1)&\mbox{if }\;|x_{1k}-a_i|\le x_{2k}-b_i\;\mbox{ and }\;
x_{2k}>b_i+r_i,\\\\
(0,-1)&\mbox{if }\;|x_{1k}-a_i|\le b_i-x_{2k}\;\mbox{ and }\;
x_{2k}<b_i-r_i,\\\\
(0,0)&\mbox{otherwise.}
\end{array}\right.
\end{eqnarray*}}
\end{Corollary}
{\bf Proof}. It follows from Proposition~\ref{mk2b}, comparison between the right-hand side of \eqref{a1a} and formula \eqref{mk3}, and the square calculations of Corollary~\ref{tosquares}. $\h$\vspace*{0.05in}

The following two examples present implementations of the subgradient algorithm realization from Corollary~\ref{mk2b} in the generalized Heron problem under consideration with
square and ball constraints, respectively.
{\small\begin{figure}[h]
\vspace{-30pt}
\begin{minipage}{2in}
   \vspace{10pt}
   \includegraphics[width=4.3in]{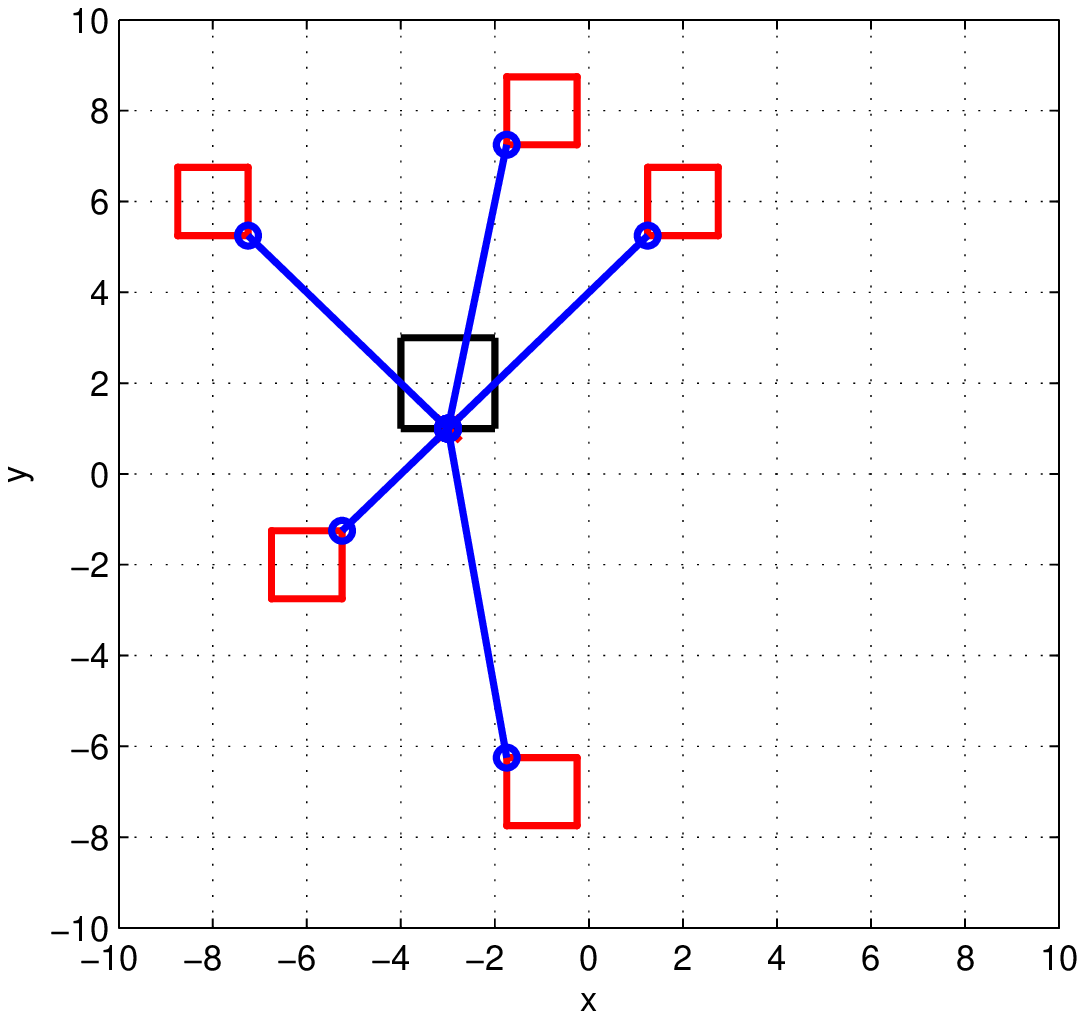}\\
\end{minipage}
~\hfill~
\begin{minipage}[t]{0.52\textwidth}
\begin{tabular}{|c|c|c|}
\hline
\multicolumn{3}{|c|}{MATLAB RESULT} \\
\hline
$k$ & $x_k$ & $V_k$ \\
\hline
1      & (-4,3)             & 26.25000 \\
10     & (-3.12500,1.04603) & 24.37500 \\
100    & (-2.99136,1.00070) & 24.25068 \\
1,000  & (-3.00133,1.00133) & 24.25000 \\
10,000 & (-2.99996,1.00001) & 24.25000 \\
15,000 & (-3.00013,1.00007) & 24.25000 \\
20,000 & (-3.00000,1.00000) & 24.25000 \\
25,000 & (-3.00001,1.00001) & 24.25000 \\
30,000 & (-3.00001,1.00001) & 24.25000 \\
\hline
\end{tabular}
\end{minipage}
\vspace{-25pt}
\caption{A Generalized Heron Problem for Squares with a Ball
Constraint with Respect to ``Max" Distances.}
\end{figure}}

{\small\begin{figure}[h]
\vspace{-30pt}
\begin{minipage}{2in}
   \vspace{10pt}
   \hspace{10pt}\includegraphics[width=4.2in]{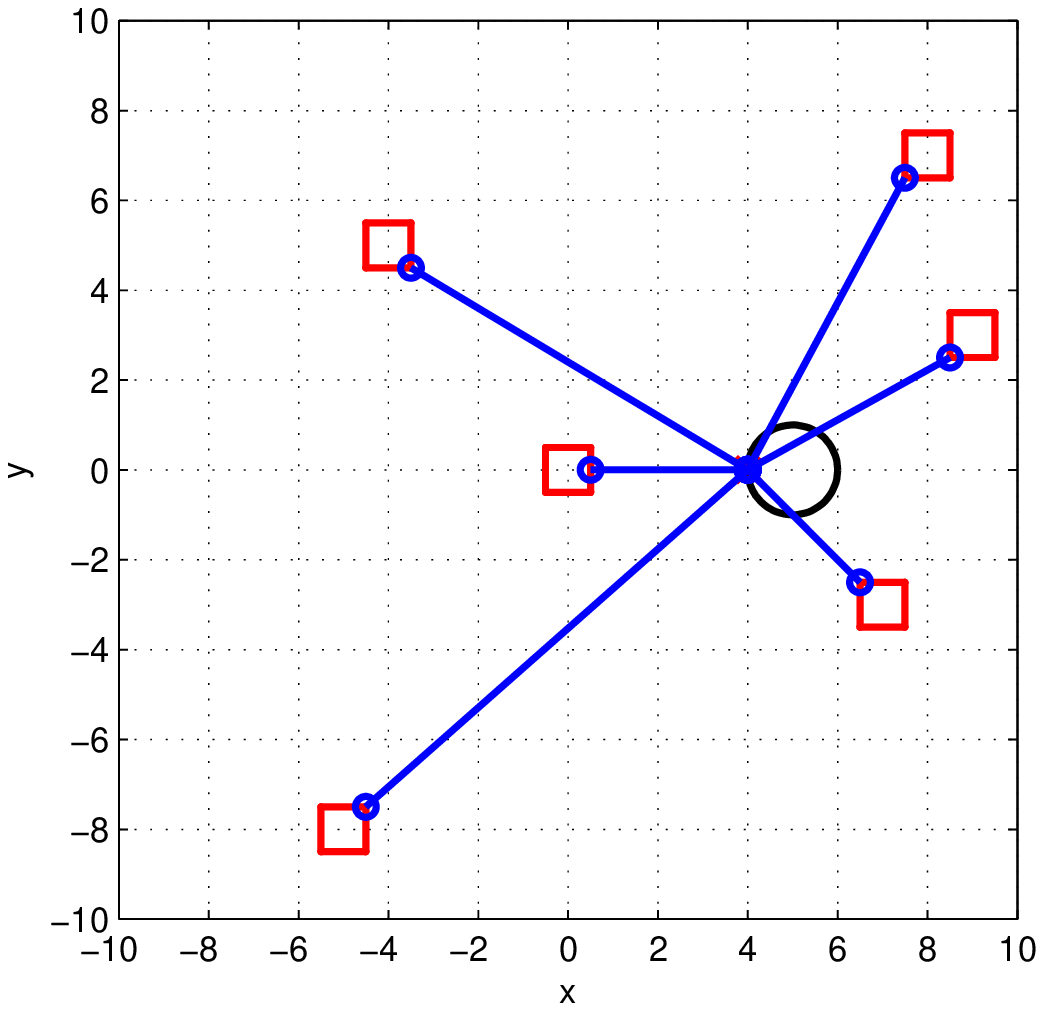}\\
\end{minipage}
~\hfill~
\begin{minipage}[t]{0.52\textwidth}
\begin{tabular}{|c|c|c|}
\hline
\multicolumn{3}{|c|}{MATLAB RESULT} \\
\hline
$k$ & $x_k$ & $V_k$ \\
\hline
1     & (5,0)             & 35 \\
10    & (4.00062,0.03519) & 33 \\
100   & (4.00000,0.00038) & 33 \\
200   & (4.00000,0.00010) & 33 \\
400   & (4.00000,0.00002) & 33 \\
600   & (4.00000,0.00001) & 33 \\
800   & (4.00000,0.00001) & 33 \\
1,000 & (4.00000,0.00000) & 33 \\
1,200 & (4.00000,0.00000) & 33 \\
\hline
\end{tabular}
\end{minipage}
\vspace{-20pt}
\caption{A Generalized Heron Problem for Squares with Ball
Constraint with Respect to ``Max" Distances.}
\end{figure}}
\begin{Example}\label{ex4.7} {\bf (generalized Heron problem with square dynamics, targets, and constraints).}
{\rm Consider the implementation of the algorithm from Corollary~{\rm\ref{mk2b}} in problem \eqref{ft} with
the square constraint $\Omega$ of center (-3,2) and short radius 1 and
the target square sets $\Omega_i$ as $i=1,\ldots,5$ of centers (-8,6), (-6,-2), (-1,8),
(-1,-7), and (2,6) with the same short radius $r=0.75$. In Figure~7 we present the results of calculations by \eqref{al} with
$\al_k=1/k$ and the starting point $x_1= (-4,3)$. The optimal solution here is $\ox\approx(-3.00001,1.00001)$ and the
optimal value is $\Hat V\approx24.25000$.}
\end{Example}
\begin{Example}\label{ex4.8} {\bf (generalized Heron problem with square dynamics and targets and with ball constraints).}
{\rm Consider the implementation of the subgradient algorithm from Corollary~\ref{mk2b} in problem \eqref{ft} with the square dynamics, the square targets $\Omega_i$ as $i=1,\ldots,6$
of centers (-5,-8), (-4,5), (0,0), (8,7), (9,3), and (7,-3) with the same short radius $r=0.5$, and with the ball constraint $\Omega$ of center (5,0) and
radius 1. The presented calculations are performed by \eqref{al} with
$\al_k=1/k$ and the starting point $x_1=(5,0)$; see Figure~8. The obtained optimal solution is $\ox\approx(4.00000,0.00000)$ with the
optimal value $\Hat V\approx 33.00000$.}
\end{Example}

Our last example concerns a three-dimensional distance version of the generalized Heron problem \eqref{distance function} for {\em cubes} of right position in $\R^3$ subject to a
ball constraint.
{\small\begin{figure}[here]
\vspace{-10pt}
\begin{minipage}{2in}
   \includegraphics[width=4in]{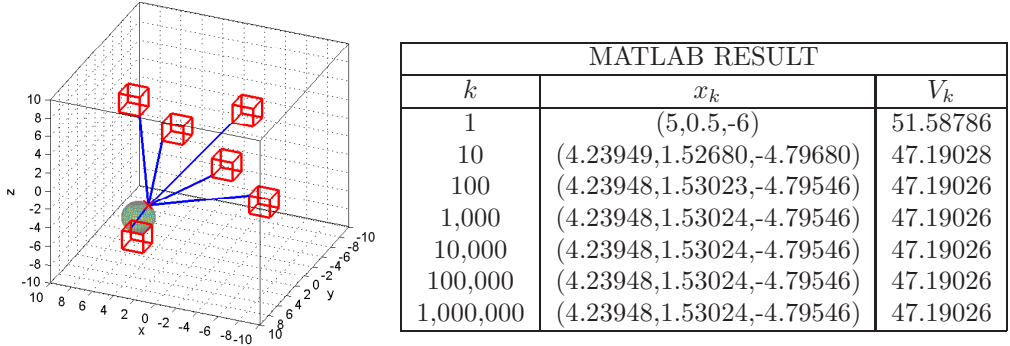}\\
\end{minipage}
~\hfill~
\begin{minipage}[t]{0.6\textwidth}
\begin{tabular}{|c|c|c|}
\hline
\multicolumn{3}{|c|}{MATLAB RESULT} \\
\hline
$k$ & $x_k$ & $V_k$ \\
\hline
1         & (5,0.5,-6)                 & 51.58786 \\
10        & (4.23949,1.52680,-4.79680) & 47.19028 \\
100       & (4.23948,1.53023,-4.79546) & 47.19026 \\
1,000     & (4.23948,1.53024,-4.79546) & 47.19026 \\
10,000    & (4.23948,1.53024,-4.79546) & 47.19026 \\
100,000   & (4.23948,1.53024,-4.79546) & 47.19026 \\
1,000,000 & (4.23948,1.53024,-4.79546) & 47.19026 \\
\hline
\end{tabular}
\end{minipage}
\vspace{-20pt}
\caption{A Generalized Heron Problem for Cubes with Ball Constraint
in Three Dimensions.}
\end{figure}}
\begin{Example}\label{ex4.9} {\bf (generalized Heron problem for cubes with ball constraints).}
{\rm Consider problem \eqref{distance function} for cubes $\Omega_i$ as $i=1,\ldots,6$ of right position in $\R^3$ with the centers $(8,-4,3)$, $(-2,-6,3)$,
$(3,-2,2)$, $(-4,-5,-6)$, $(-3,1,1)$, and $(3,7,-5)$  and the same
short radius 1 subject to the ball constraint $\Omega$ of center (5,2,-6) and radius 1.5. The projection $P((x,y,z);\Omega)$ and quantities $q_{ik}$ in algorithm \eqref{al} are calculated
similarly to Example~\ref{ex4.1}. Figure~9 presents the implementation of the subgradient algorithm \eqref{al} with $\al_k =1/k$  and the starting point $x_1=(5,5,-6)$.
As we see, the optimal solution calculated here up to five
significant digits is $\ox\approx(4.23948,1.53024,-4.79546)$ and the optimal value is
$\Hat V\approx47.19026$.}
\end{Example}

We conclude the paper by the following three observations.

\begin{Remark}\label{rem} {\bf (extensions and other location problems).}
{\rm

{\bf (i)} Note that the approach and results of this paper can be easily extended to the {\em weighted
version} of the generalized Heron problem \eqref{ft}:
\begin{equation}\label{ft1}
\mbox{minimize }\;T(x):=\sum_{i=1}^n\mu_i
T^F_{\Omega_i}(x),\;\mbox{subject to }\;x\in\Omega,
\end{equation}
where $\mu_i\ge 0$ as $i=1,\ldots,n$ are given weights. Since we have
\begin{equation*}
\partial\big(\mu_i T^F_{\Omega_i}\big)(\ox)=\mu_i \partial T^F_{\Omega_i}(\ox)
\end{equation*}
for both convex and nonconvex subdifferentials used in this paper, it is straightforward to derive counterparts of the qualitative and numerical results obtained above for the case of the weighted generalized Heron problem \eqref{ft1}. For example, the equation
\begin{equation*}
\sum_{i=1}^n\mu_i\cos\big(a_i(\ox),v\big)=0\;\mbox{ for every }\;v\in L(\ox)\setminus\{0\}
\end{equation*}
replaces the one in (\ref{cos rep}) for all the corresponding results.

{\bf (ii)} Our variational approach can be used to solve a variety of other
facility location problems. In particular, the following {\em smallest intersecting ball problem} can be naturally formulated and investigated by using the above tools of variational analysis and generalized differentiation: given $n$ nonempty closed
subsets $\Omega_i\subset X,\;i=1,\ldots,n$, find a point $\ox$ on a given set $\Omega$
and the smallest number $r>0$ such that the ball with center at
$\ox$ and radius $r$ has nonempty intersection with all the sets $\Omega_i$ as
$i=1,\ldots,n$. This problem is modeled as follows:
\begin{equation*}
\mbox{minimize
}\;\mbox{M}(x):=\max\big\{d(x;\Omega_i)\big|\;i=1,\ldots,n\big\}\;
\mbox{ subject to }\;x\in \Omega.
\end{equation*}
We intend to address this and other facility location problems in our future research.

{\bf (iii)} For some results in the Hilbert space setting of Section 3, it is possible to use \emph{the proximal normal cone} instead of the Fr\'echet normal cone. However, we use the Fr\'echet normal cone consistently for the simplicity of presentation. }
\end{Remark}

{\bf Acknowledgments.} The authors are grateful to Jon Borwein, Mari\'an Fabian, and Doan The Hieu for valuable discussions on the material of this paper. The first author
acknowledges partial supports by the USA National Science Foundation under grant
DMS-1007132, by the European Regional Development Fund (FEDER), and by the following Portuguese agencies: Foundation for Science and Technologies (FCT),
Operational Program for Competitiveness Factors (COMPETE), and Strategic Reference Framework (QREN).

\end{document}